\documentclass[12pt,a4paper]{amsart}
\usepackage[utf8]{inputenc}
\usepackage[T1]{fontenc}
\usepackage{microtype}
\usepackage{amsmath,amsfonts,amssymb,amsthm}
\numberwithin{equation}{section}
\usepackage{ifthen,xifthen}
\usepackage[shortlabels]{enumitem}
\usepackage[centering]{geometry}
\usepackage{hyperref}
\hypersetup{colorlinks,allcolors=black,final}
\usepackage[capitalize,noabbrev,nameinlink]{cleveref}
\usepackage[shortlabels]{enumitem}
\usepackage{tikz}
\usetikzlibrary{external}
\write18{if [ ! -d figures ]; then mkdir figures; fi}
\tikzsetexternalprefix{figures/}

\usepackage{ragged2e,longtable,tabu,booktabs}

\newtheorem{Theorem}{Theorem}[section]\Crefname{Theorem}{Theorem}{Theorems}
\Crefname{Lemma}{Lemma}{Lemmas}
\newtheorem{Proposition}[Theorem]{Proposition}\Crefname{Proposition}{Proposition}{Propositions}
\Crefname{Corollary}{Corollary}{Corollaries}
\theoremstyle{definition}
\newtheorem{Definition}[Theorem]{Definition}\Crefname{Definition}{Definition}{Definitions}
\Crefname{Remark}{Remark}{Remarks}
\newtheorem{Example}[Theorem]{Example}\Crefname{Example}{Example}{Examples}

\providecommand{\abs}[1]{\lvert#1\rvert}

\providecommand{\Bigabs}[1]{\Bigl|#1\Bigr|}
\providecommand{\biggabs}[1]{\biggl|#1\biggr|}

\providecommand{\norm}[2][]{\lVert#2\rVert\ifthenelse{\isempty{#1}}{}{_{#1}}}
\providecommand{\bignorm}[2][]{\bigl\|#2\bigr\|\ifthenelse{\isempty{#1}}{}{_{#1}}}
\providecommand{\Bignorm}[2][]{\Bigl\|#2\Bigr\|\ifthenelse{\isempty{#1}}{}{_{#1}}}
\providecommand{\biggnorm}[2][]{\biggl\|#2\biggr\|\ifthenelse{\isempty{#1}}{}{_{#1}}}
\providecommand{\Biggnorm}[2][]{\Biggl\|#2\Biggr\|\ifthenelse{\isempty{#1}}{}{_{#1}}}

\providecommand{\TVnorm}[1]{\norm[\text{TV}]{#1}}
\providecommand{\ceil}[1]{\lceil#1\rceil}

\providecommand{\biggceil}[1]{\biggl\lceil#1\biggr\rceil}

\providecommand{\floor}[1]{\lfloor#1\rfloor}

\providecommand{\cA}{\mathcal A}
\providecommand{\cB}{\mathcal B}
\providecommand{\Borel}{\cB}
\providecommand{\R}{\mathbb R}
\providecommand{\Z}{\mathbb Z}
\providecommand{\ifu}[1]{\mathbf 1_{#1}}
\providecommand{\Ioc}[2]{(#1,#2]}
\providecommand{\Ioo}[2]{(#1,#2)}
\providecommand{\Ico}[2]{[#1,#2)}
\providecommand{\Icc}[2]{[#1,#2]}
\providecommand{\spr}[2]{\langle#1,#2\rangle}
\providecommand{\xto}{\xrightarrow}
\let\originald\d
\renewcommand{\d}{\ifthenelse{\boolean{mmode}}{\mathrm d}{\originald}}
\providecommand{\dd}{\,\d}

\providecommand{\from}{\colon}
\providecommand{\argmt}{{\mathchoice{{}\cdot{}}
                                    {{}\cdot{}}
                                    {{}\bullet{}}
                                    {{}\bullet{}}}}
\providecommand{\Laplace}{\Delta}
\providecommand{\textq}[1]{\text{#1}\quad}
\providecommand{\qtext}{\quad\text}
\providecommand{\qtextq}[1]{\quad\text{#1}\quad}
\providecommand{\Prob}{\mathbb P}
\providecommand{\E}{\mathbb E}
\providecommand{\tensor}{\otimes}
\providecommand{\Tensor}{\bigotimes}
\providecommand{\cF}{\mathcal F}
\DeclareMathOperator{\supp}{supp}
\DeclareMathOperator{\tr}{tr}
\DeclareMathOperator{\dist}{dist}
\DeclareMathOperator{\diam}{diam}
\providecommand{\setsize}{\abs}
\providecommand{\N}{\mathbb N}
\providecommand{\B}{\mathbb B}

\providecommand{\isect}{\cap}
\providecommand{\Union}{\bigcup}
\providecommand{\union}{\cup}
\providecommand{\cV}{\mathcal V}
\providecommand{\cE}{\mathcal E}
\providecommand{\Hvis}{H_{\textrm{vis}}}

\providecommand{\Psite}{\Prob_{\mathit s}}
\providecommand{\Pedge}{\Prob_{\mathit e}}
\providecommand{\floor}[1]{\lfloor#1\rfloor}
\providecommand{\Fbar}{\overline F}
\providecommand{\fbar}{\overline f}
\providecommand{\ve}{\varepsilon}

\DeclareMathOperator{\id}{id}

\providecommand{\constb}[1][f]{D_{#1}}
\providecommand{\Normal}{\mathcal N}
\providecommand{\cM}{\mathcal M}
\providecommand{\cU}{\mathcal U}

\title{Uniform existence of the IDS on lattices and groups}
\author{C.~Schumacher\and F.~Schwarzenberger\and I.~Veseli\'c}

\begin{document}

\begin{abstract}
We present a general framework for thermodynamic limits
and its applications to a variety of models.
In particular we will identify criteria such that the limits are uniform
in a parameter.
All results are illustrated with the example of eigenvalue counting functions
converging to the integrated density of states.
In this case, the convergence is uniform in the energy.
\end{abstract}

\maketitle

\tableofcontents

\section{Introduction}
\label{sec:intro}

The thermodynamic limit, i.\,e.\ taking averages over larger and larger volumes,
performs the transition from microscopic to macroscopic models.
In the context of this paper, thermodynamic limits are used to define
the integrated density of states (IDS) of discrete Schr\"odinger operators,
in particular, random ones.
The IDS of random Schr\"odinger operators has been studied in the mathematical literature at least since 1971,
see e.\,g.\ \cite{Pastur-71,Shubin-79,KirschM-07,Veselic-08}.
The IDS measures the number of quantum states per unit volume
below a threshold energy~$E$ and is thus a function of~$E$.
The rigorous implementation of this notion involves a thermodynamic limit,
see \cref{sec:defsandmodels}.
Classical results about the existence of this macroscopic volume limit
apply to fixed energy, which is to say that the limit is considered pointwise.
More recently, there has been increased interest to study 
stronger forms of convergence, for instance uniformly in the energy parameter.
This has been studied in various geometric and stochastic contexts.
For instance, there are works devoted to Hamiltonians associated to quasicrystals, 
\cite{LenzS-03a,LenzS-06}, to graphings \cite{Elek-07} and percolations graphs
\cite{Veselic-05b,Veselic-06,AntunovicV-08b-short,SamavatSV-14},
as well as abstracts frameworks covering general classes of examples 
\cite{LenzV-09}, where this is certainly a non-exhaustive list of references. 
We will review here in particular the results obtained in 
\cite{LenzMV-08},
\cite{LenzSV-10,LenzSV-12}, 
\cite{PogorzelskiS-16}, 
\cite{SchumacherSV-17}, and 
\cite{SchumacherSV-18}. They have all in common that their findings can be formulated as Banach-valued 
ergodic theorems. 

The convergence results in this work are not restriced to the IDS of random Schr\"odinger operators,
but provide a general framework for thermodynamic limits with respect to a Banach space topology.
To this end, we introduce almost additive fields,
which, after normalization with volume,
converge as the volume exhausts the physical space.
The microscopic structure, like the values of the potential of the
Schr\"odinger operator and/or whether a percolation site is open or closed,
are encoded in a \emph{coloring}.
In applications, the coloring is often random,
for the lack of detailed knowledge about the microscopic properties of the material.

For the thermodynamic limit to exist, one needs a certain homogeneity of the coloring.
In our case, finite portions of the coloring, called \emph{patterns},
should occur with a certain frequency.
If there are only finitely many \emph{colors}, i.\,e.\ values of the coloring,
this condition is natural and allows the existence of the thermodynamic limit
to be proven in any Banach space.
For random Schr\"odinger operators, one chooses the Banach space
of right continuous bounded functions with $\infty$-norm.
Unfortunately, typical potentials have infinitely many values.
In that setting, we model the colorings as random fields with some independence.
This allows to employ a multivariate version of the theorem of Glivenko--Cantelli.

The paper is structured as follows.
In \cref{sec:defsandmodels}, we introduce some specific random
Schr\"odinger operators and their IDS in more detail.
They serve as examples for the abstract theorems
as well as illustrations of the meaning of our assumptions.
We then present a series of theorems on thermodynamic limits
with increasing complexity.
\Cref{sec:finiteA} contains the oldest of the presented results, which was in some sense the motivation for further developments.
It considers fields over~$\Z^d$ obeying the finiteness condition~\eqref{eq:Afinite}
on the set of colors.
Since the results are easier to motivate and less technical than later ones, they serve well as a point of departure.

In \cref{sec:amenable},
we introduce briefly finitely generated amenable groups
and show how the main result generalizes to this setting.
We first consider amenable groups which satisfy an additional tiling condition.
To remove this tiling condition, we outline quasi tilings
and discuss the result on amenable groups,
still keeping the finiteness condition.

The finiteness condition~\eqref{eq:Afinite}
is finally tackled and removed in \cref{sec:infiniteA}.
For clarity, we first deal with the euclidean lattice~$\Z^d$
and only after that exhibit the most general formulation on
finitely generated amenable groups without finiteness condition.

An outlook, a list of symbols, and a list of references conclude the paper.

\section{Physical models}
\label{sec:defsandmodels}
In this section, we introduce example systems to motivate the abstract results.
For the moment we will stick to the simple geometry of~$\Z^d$
but use a notation that later generalizes
to general geometries on amenable groups in a straightforward way.

\subsection{The Anderson model on \texorpdfstring{$\Z^d$}{Zd}}
\label{sec:Anderson_Zd}

The Anderson model, introduced by P. W. Anderson in \cite{Anderson-58},
is a prototypical random Schr\"odinger operator.
To define it, let us introduce some notation.

As physical space we choose the group $G:=\Z^d$. 
To make neighborhood relations explicit, we introduce the
\emph{Cayley graph} of~$G$.
The Cayley graph of~$G$ has~$G$ itself as vertex set,
and two vertices $v,w\in G$ are connected by an undirected edge,
if $\norm[1]{v-w}=1$, where $\norm[1]\argmt$
is the usual norm on $\Z^d$ when viewed as subset of $\ell^1\{1,\dotsc,d\}$.
In this case, we write $v\sim w$.
Thus, the considered geometry is nothing but the usual $d$-dimensional lattice.
However, to stay consistent with the notation required in later parts of the article, we already use the notion of a Cayley graph.
The general definition of Cayley graphs is introduced in \cref{sec:amenable}.

The group~$G$ acts transitively on its Cayley graph by translation:
\begin{equation*}
  \tilde\tau\colon G\times G\to G\textq,
  \tilde\tau_gv:=(g,v):=v-g\text.
\end{equation*}
This group action lifts to a unitary group action
on the square summable functions on the vertices of the Cayley graph,
$\ell^2(G)$, which we denote by
\begin{equation*}
  U_g\from\ell^2(G)\to\ell^2(G)\textq,
  (U_g\varphi)(v):=\varphi(\tilde\tau_gv)\text.
\end{equation*}

The \emph{Laplace operator} $\Laplace\from\ell^2(G)\to\ell^2(G)$, given by
\begin{equation}\label{eq:laplaceop}
  (\Laplace\varphi)(v)
  :=\sum_{w\sim v}\bigl(\varphi(w)-\varphi(v)\bigr)\text,
\end{equation}
mimics the sum of the second derivatives.
The operator~$-\Laplace$ is bounded, self-adjoint, positive semi-definite
and serves as the quantum mechanical observable of the kinetic energy.
Its spectrum is $\sigma(-\Laplace)=\Icc0{2d}$,
as one can see with Fourier analysis.
Note that the Laplace operator is equivariant with respect to $G$,
i.\,e., for all $g\in G$, we have
$\Laplace\circ U_g=U_g\circ\Laplace$.

To build a Schr\"odinger operator,
we need a multiplication operator~$V$ on $\ell^2(G)$,
which plays the role of the observable of potential energy.
At this stage, the randomness enters.
To this end, we choose a (measurable) set~$\cA\subseteq\R$,
which we call the set of colors,
equip it with the trace topology inherited from~$\R$
and its Borel $\sigma$-algebra~$\Borel(\cA)$,
and choose a probability measure~$\Prob_0$ on $(\cA,\Borel(\cA))$,
i.\,e.\ the colors.
The probability space for the random potential is
$(\Omega,\cB,\Prob)
 :=(\cA,\Borel(\cA),\Prob_0)^{\tensor G}
 :=(\cA^G,\Borel(\cA)^{\tensor G},\Prob_0^{\tensor G})$.
The \emph{random potential} $V\from\Omega\times G\to\cA$
returns for~$v\in G$ the $v$-th coordinate of~$\omega\in\Omega$:
\begin{equation}\label{eq:randompotential}
  V_\omega(v):=V(\omega,v):=\omega_v:=\omega(v)\text,
\end{equation}
so that the random variables $\Omega\ni\omega\mapsto V_\omega(v)$, $v\in G$,
are independent and identically distributed.
For each $\omega\in\Omega$, the random potential~$V_\omega$
operates on $\ell^2(G)$ by multiplication.
If~$\cA\in\Borel(\R)$ is bounded,
the multiplication operator~$V_\omega$ is bounded and self-adjoint.
Analogous to above, the group action of~$G$ on the Cayley graph
lifts to an ergodic group action
$\tau_g\from\Omega\to\Omega$ on~$\Omega$:
\begin{equation*}
  (\tau_g\omega)_v:=\omega_{\tilde\tau_gv}\text,
\end{equation*}
and the random potential is equivariant, meaning
$V_{\tau_g\omega}\circ U_g=U_g\circ V_\omega$ for all $g\in G$.

For each $\omega\in\Omega$, the Schr\"odinger operator
\begin{equation} \label{eq:schroedinger}
  H_\omega:=-\Laplace+V_\omega\from\ell^2(G)\to\ell^2(G)
\end{equation}
is well-defined, bounded, and self-adjoint.
The operator family $(H_\omega)_{\omega\in\Omega}$
is the famous \emph{Anderson Hamiltonian} and is equivariant, too:
$H_{\tau_g\omega}\circ U_g=U_g\circ H_\omega$
for all $g\in G$.
This equivariance shows that the spectrum $\sigma(H_\omega)$
is a shift invariant quantity.
Since the group action of~$G$ on~$\Omega$ is ergodic,
the spectrum of $H_\omega$ is almost surely constant:
there is a closed set $\Sigma\subseteq\R$ such that
\begin{equation*}
  \Sigma=\sigma(H_\omega)
\end{equation*}
for $\Prob$-almost all $\omega\in\Omega$.
In quantum mechanics, the spectrum is interpreted
as the set of possible energies of the particle described by~$H_\omega$.
The fact that it is deterministic
and does not depend on the microscopic structure of the material
makes the interpretation as a homogeneous material possible.
More details can be found for example in \cite{FigotinP-92,Kirsch-08}.

The main example to motivate and to illustrate the results
in later chapters is the integrated density of states (IDS),
also known as the spectral distribution function,
of the Anderson Hamiltonian.
Its definition is
\begin{equation*}
  N\from\R\to\R\textq,
  N(E):=\E[\spr{\delta_0}{\ifu{\Ioc{-\infty}E}(H_\omega)\delta_0}]\text,
\end{equation*}
where $\delta_v\in\ell^2(G)$ is the Kronecker delta at $v\in G$,
i.\,e.\ $\delta_0=\ifu{\{v\}}$,
and $\ifu{\Ioc{-\infty}E}(H_\omega)$
is the spectral projection of~$H_\omega$
onto the energy interval $\Ioc{-\infty}E$.
The IDS is monotone, bounded by~$1$,
and gives rise to a probability measure~$\d N$,
called the \emph{density of states measure}.
The topological support of the density of states
is the almost sure spectrum of~$H_\omega$.
In the following,
we will describe how to approximate the IDS using only finite matrices.

For each $\Lambda\subseteq G$, we write $\ell^2(\Lambda)$
for the subspace of $\varphi\in\ell^2(G)$
with support $\supp\varphi\subseteq\Lambda$.
The indicator function of~$\Lambda$, used as multiplication operator on~$\ell^2(G)$,
is the self-adjoint orthogonal projection $\ifu\Lambda\from\ell^2(G)\to\ell^2(\Lambda)$.
The operator $H_\omega^\Lambda\from\ell^2(\Lambda)\to\ell^2(\Lambda)$ is given by
\begin{equation*}
  H_\omega^\Lambda:=\ifu\Lambda\circ H_\omega\circ(\ifu\Lambda)^*\text.
\end{equation*}

Now let~$\cF$ be the (countable) set of finite subsets of~$G$
and assume $\Lambda\in\cF$.
The representing matrix of~$H_\omega^\Lambda$
contains the matrix elements $\spr{\delta_v}{H_\omega\delta_w}$
with $v,w\in\Lambda$ and is thus the $\Lambda\times\Lambda$ clipping of the representing
$G\times G$ matrix of~$H_\omega$.
Of course, $H_\omega^\Lambda$ is Hermitian and has real eigenvalues.
For $\omega\in\Omega$, $\Lambda\in\cF$,
the eigenvalue counting function $n(\Lambda,\omega)\from\R\to\R$,
\begin{equation*}
  n(\Lambda,\omega)(E):=\tr\ifu{\Ioc{-\infty}E}(H_\omega^\Lambda)
    =\dim\{\varphi\in\ell^2(\Lambda)\mid\spr\varphi{H_\omega^\Lambda\varphi}\le E\}
\end{equation*}
is continuous from the right and counts the eigenvalues of~$H_\omega^\Lambda$
below the threshold energy~$E$ according to their multiplicity.
We will view the family of eigenvalue counting functions as
\begin{equation*}
  n\from\cF\times\Omega\to\B\text,
\end{equation*}
where~$\B$ is the Banach space of right continuous functions
from~$\R$ to~$\R$.

Denote by $\Lambda_L:=\Ico0L^d\isect\Z^d\in\cF$ a cube of side length $L\in\N$.
The celebrated \emph{Pastur--Shubin formula} states that,
for $\Prob$-almost all $\omega\in\Omega$
and all $E\in\R$ where the IDS~$N$ is continuous,
we have
\begin{equation}\label{eq:PasturShubin}
  N(E)
  =\lim_{L\to\infty}\setsize{\Lambda_L}^{-1}n(\Lambda_L,\omega)(E)\text,
\end{equation}
see e.\,g.\ \cite{SchumacherS-15}.
In fact, since for this model in particular the IDS~$N$
is continuous at all energies, a straight forward argument,
using also the boundedness and the monotonicity of~$N$,
shows that the convergence does not only hold for all~$E\in\R$,
but is actually uniform in~$E$.
The Pastur--Shubin formula can be viewed as an ergodic theorem,
since it states the equality of an ensemble average and a spatial average.

\subsection{Quantum percolation models}
\label{sec:percolation}

\subsubsection{Site and edge percolation}
\label{sec:siteedge}

We first introduce \emph{site percolation} on the Cayley graph of $G=\Z^d$.
Let $\cA:=\{0,1\}$, and as in \cref{sec:Anderson_Zd} let
$(\Omega,\cB,\Prob):=(\cA,\Borel(\cA),\Psite)^{\tensor G}$
with a (non-de\-ge\-ne\-rate) Bernoulli probability~$\Psite$ on~$\cA$.
In the percolation setting, the randomness does not determine the potential,
but configuration space itself.
The \emph{site percolation graph} $(\cV_\omega,\cE_\omega)$
corresponding to $\omega\in\Omega$
is induced by the Cayley graph of~$G$ on the vertex set
\begin{equation*}
  \cV_\omega:=\{v\in G\mid\omega_v=1\}\text.
\end{equation*}
This means that the edge set is
\begin{equation*}
  \cE_\omega:=\{\{v,w\}\subseteq\cV_\omega\mid v\sim w\}\text,
\end{equation*}
that is, we keep the edges of the Cayley graph
which have both end points in~$\cV_\omega$.

In \emph{edge percolation}, one does not erase random sites
from the Cayley graph but instead random edges.
The color set $\cA:=\{0,1\}^d$ is suitable to implement this strategy.
Namely, for a (non-degenerate) Bernoulli probability~$\Pedge$ on $\{0,1\}$,
we define
$(\Omega,\cB,\Prob):=(\cA,\Borel(\cA),\Pedge^{\tensor d})^{\tensor G}$
and define the edge set corresponding to $\omega\in\Omega$ by
\begin{equation*}
  \cE_\omega:=\Union\nolimits_{j=1}^d\bigl\{\{v,v+e_j\}\bigm|(\omega_v)_j=1\bigr\}
\end{equation*}
where~$e_j$ is the $j$-th standard basis vector of~$\Z^d$.
The vertex set of the edge percolation graph contains all vertices
to which an edge in~$\cE_\omega$ is attached:
\begin{equation}\label{eq:verticesInEdgePercolation}
  \cV_\omega:=\{v\in G\mid\exists w\in G\colon\{v,w\}\in\cE_\omega\}\text.
\end{equation}

The Laplace operator on a subgraph $(\cV_\omega,\cE_\omega)$
of the Cayley graph of~$G$ is
\begin{equation}\label{eq:laplace}
  \Laplace_\omega\from\ell^2(\cV_\omega)\to\ell^2(\cV_\omega)\textq,
  (\Laplace_\omega\varphi)(v):=
  \sum_{w\in\cV_\omega,\{v,w\}\in\cE_\omega}\bigl(\varphi(w)-\varphi(v)\bigr)
  \text.
\end{equation}
Since we want to use the group action of~$G$,
we define the Hamiltonians on $\ell^2(G)$ via
\begin{equation*}
  H_\omega\from\ell^2(G)\to\ell^2(G)\textq,
  (H_\omega\varphi)(v):=
  \begin{cases}
    -(\Laplace_\omega\varphi)(v)&\text{if $v\in\cV_\omega$ and}\\
    \alpha\varphi(v)            &\text{if $v\in G\setminus\cV_\omega$}
  \end{cases}
\end{equation*}
with a constant $\alpha\in\Ioo{2d}\infty$.
The operator~$H_\omega$ leaves the subspaces $\ell^2(\cV_\omega)$
and $\ell^2(G\setminus\cV_\omega)$ invariant.
Since $\Laplace_\omega$ is bounded by~$2d$ and $\alpha>2d$,
the spectrum of of~$H_\omega$ is the disjoint union of
$\sigma(-\Laplace_\omega)$ and $\{\alpha\}$,
and the two components can be studied separately.

In fact, the percolation Hamiltonians are equivariant,
and again their spectrum is almost surely constant.
As in \cref{sec:Anderson_Zd}, we define the IDS $N\from\R\to\R$
and the eigenvalue counting function $n\from\cF\times\Omega\to\B$
of~$H_\omega$.
The Pastur--Shubin formula~\eqref{eq:PasturShubin}
remains correct $\Prob$-almost surely for all energies $E\in\R$
at which~$N$ is continuous, see for instance
\cite{Veselic-05a,Veselic-05b,KirschM-06,Veselic-06,AntunovicV-08b-short}.
But in contrast to the Anderson Hamiltonian,
the IDS of a percolation Hamiltonian is not continuous.
In fact, the set of discontinuities of~$N$ is dense in~$\Sigma$.
Note that since the IDS is monotone,
so there can be at most countably many discontinuities.

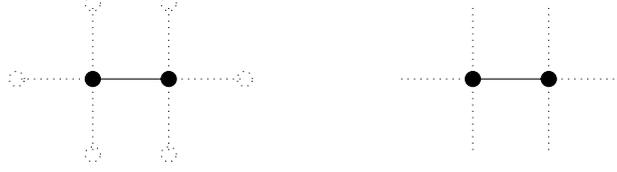
\begin{figure}
\begin{center}
\tikzsetnextfilename{prob4cluster}%
\begin{tikzpicture}
  \draw[fill] (0,0) circle[radius=1mm] -- (1,0) circle[radius=1mm];
  \foreach \x/\y in {0/1, -1/0, 0/-1}
    \draw[dotted] (0,0) -- (\x,\y) circle[radius=1mm];
  \foreach \x/\y in {1/1, 2/0, 1/-1}
    \draw[dotted] (1,0) -- (\x,\y) circle[radius=1mm];
  \begin{scope}[xshift=5cm]
  \draw[fill] (0,0) circle[radius=1mm] -- (1,0) circle[radius=1mm];
    \foreach \x/\y in {0/1, -1/0, 0/-1}
      \draw[dotted] (0,0) -- (\x,\y);
    \foreach \x/\y in {1/1, 2/0, 1/-1}
      \draw[dotted] (1,0) -- (\x,\y);
  \end{scope}
\end{tikzpicture}
\end{center}
\caption{The probability that $0\in G$ is contained in a cluster of size~$2$
  is $2d\Psite(1)^2\Psite(0)^{4d-2}$ in site percolation
  and $2d\Pedge(1)\Pedge(0)^{4d-2}$ in edge percolation.
  Illustration for $d=2$.}
\label{fig:clustersize2}
\end{figure}%
That the IDS is discontinuous, can be seen as follows.
The percolation graph splits into its connected components,
the so-called \emph{clusters}.
More precisely, the clusters of a graph are the equivalence classes
of the minimal equivalence relation for which neighbors are equivalent.
The event that $0\in G$
is contained in a finite cluster has positive probability,
see \cref{fig:clustersize2}.
A finite cluster of size~$s\in\N$ supports~$s$ eigenfunctions of~$H_\omega$.
In particular, the constant function with value $s^{-1/2}$
is an $\ell^2$-normalized eigenfunction of~$H_\omega$ with eigenvalue~$0$.
On the event that~$0$ is contained in a finite cluster of size~$s$,
we have $\spr{\delta_0}{\ifu{\Ioc{-\infty}0}(H_\omega)\delta_0}=s^{-1}$.
Thus
\begin{equation*}
  N(0)
  =\E[\spr{\delta_0}{\ifu{\Ioc{-\infty}0}(H_\omega)\delta_0}]
  \ge\sum_{s\in\N}s^{-1}\Prob(\text{the cluster of~$0$ has size~$s$})
  >0\text,
\end{equation*}
and since $N(E)=0$ for all $E<0$, we verified that~$N$ is discontinuous.

Note that there can be compactly supported eigenfunctions of~$H_\omega$
on infinite clusters.
To construct an example, consider a finite symmetric cluster like
$\{-e_1,0,e_1\}\subseteq G$.
The symmetry that flips the sign of~$e_1$
commutes with the Laplace operator restricted to this cluster.
Therefore, there are antisymmetric eigenfunctions like $(-1,0,1)/\sqrt2$,
which have to vanish on each fixed point of the symmetry.
Now connect the finite symmetric cluster to an infinite cluster
with edges only touching the fixed points.
For more details see \cite{ChayesCF-85a,Veselic-05b}.


\subsubsection{The Anderson model on a percolation graph}
\label{sec:AndPerc}
We will now combine the random kinetic energy of the percolation Hamiltonian
with the random potential energy of the Anderson model.
A suitable set of colors is $\cA=\cA_0\times\{0,1\}$
for site and $\cA=\cA_0\times\{0,1\}^d$ for edge percolation,
where the bounded set~$\cA_0\subseteq\R$
contains the values of the random potential.
Accordingly, the probability space is
\begin{equation*}
  (\Omega,\cB,\Prob):=
  \begin{cases}
    (\cA,\Borel(\cA),\Prob_0\tensor\Psite)^{\tensor G}&\text{or}\\
    (\cA,\Borel(\cA),\Prob_0\tensor\Pedge^{\tensor d})^{\tensor G}\text.
  \end{cases}
\end{equation*}
For each $v\in G$, the random parameter $\omega=(\omega_v)_{v\in G}\in\Omega$
has a coordinate $\omega_v=(\omega'_v,\omega''_v)$ with $\omega'_v\in\cA_0$
and either $\omega''_v\in\{0,1\}$ or $\omega''_v\in\{0,1\}^d$.
The random potential~$V\from G\to\cA_0\subseteq\R$ is defined by
\begin{equation}
  V_\omega(v):=V(\omega,v):=\omega'_v\text.
\end{equation}
Analogously to the construction in \cref{sec:siteedge},
we define the site percolation graph as the graph induced
by the Caley graph of~$G$ on the vertex set
\begin{equation*}
  \cV_\omega:=\{v\in G\mid\omega''_v=1\}\text.
\end{equation*}
The edge percolation graph has edge set
\begin{equation*}
  \cE_\omega:=\Union\nolimits_{j=1}^d
    \bigl\{\{v,v+e_j\}\bigm|(\omega''_v)_j=1\bigr\}
\end{equation*}
and vertex set given by~\eqref{eq:verticesInEdgePercolation}.
The Laplace operator on the percolation graph is given by~\eqref{eq:laplace}.
The Hamiltonians $H_\omega\from\ell^2(G)\to\ell^2(G)$ are defined by
\begin{equation*}
  (H_\omega\varphi)(v):=
  \begin{cases}
    -(\Laplace_\omega\varphi)(v)+V_\omega(v)
      &\text{if $v\in\cV_\omega$ and}\\
    \alpha\varphi(v)
      &\text{if $v\in G\setminus\cV_\omega$}
  \end{cases}
\end{equation*}
with $\alpha>2d+\sup\cA_0$.
Again, the spectrum $\sigma(H_\omega)$ will have a component~$\{\alpha\}$
and a disjoint remainder, which is the part we are most interested in.
We have again an equivariant group action, namely for all $g\in G$,
\begin{equation*}
  H_{\tau_g\omega}\circ U_g
  =U_g\circ H_{\tau_g\omega}\text.
\end{equation*}

If the potential has a probability density,
the randomness will smooth out at least some discontinuities of the IDS,
see \cite{Veselic-06}.

\section{Ergodic theorems for finite colors on~\texorpdfstring{$\Z^d$}{Zd}}
\label{sec:finiteA}

In \cite[Theorem~2]{LenzMV-08},
the authors prove a quantitative version of the following statement.
\begin{Theorem}\label{thm:LMV-ecf}
  In either of the settings presented in \cref{sec:defsandmodels},
  assume that the set~$\cA$ of colors is finite:
  \begin{equation}\label{eq:Afinite}
    \setsize\cA<\infty\text.
  \end{equation}
  Then, 
  for $\Prob$-almost all $\omega\in\Omega$,
  \begin{equation*}
    \lim_{L\to\infty}\norm[\infty]
      {\setsize{\Lambda_L}^{-1}n(\Lambda_L,\omega)-N}
    =0\text,
  \end{equation*}
  where $N\from\R\to\R$ is the IDS.
\end{Theorem}
Condition~\eqref{eq:Afinite} is automatically satisfied,
unless a potential with infinitely many values is present.

In \cite{LenzMV-08},
the authors do not talk about a probability space with many configurations,
but rather fix one configuration and argue on some required properties.
This properties turn out to be almost surely satisfied in the examples in
\cref{sec:defsandmodels}.
Nonetheless, it is illuminating to see a deterministic example,
which we present next.

\subsection{Visible points}
\label{sec:visible}
A point of~$\Z^d$ is \emph{visible} from the origin,
if there is no other point of~$\Z^d$
on the straight line connecting the origin and the point in question,
see \cref{fig:visiblepoints_small,fig:visiblepoints_large}.
The set of visible points is thus
\begin{equation*}
  \cV:=\{x\in\Z^d\mid\{tx\mid t\in\Icc01\}\isect\Z^d=\{0,x\}\}\text.
\end{equation*}
Equivalently, a point $x\in\Z^d$ is visible, if $x=0$
or if the greatest common divisor of its coordinates is~$1$.
See \cite{BaakeMP-00} for a systematic exploration of~$\cV$.

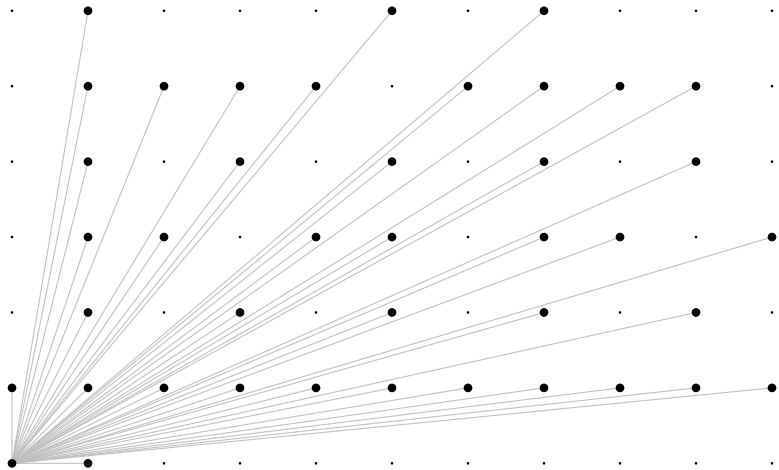
\begin{figure}
\pgfmathsetmacro{\X}{10}
\pgfmathsetmacro{\Y}{6}
\begin{center}
\tikzsetnextfilename{visible1}%
\begin{tikzpicture}
  \foreach \y in {0, ..., \Y} {
    \foreach \x in {0, ..., \X} {
      \pgfmathsetmacro{\visible}{ gcd(\x,\y) == 1 }
      \ifthenelse{\visible = 1}
        {\draw[draw=gray!50] (0,0) -- (\x,\y);
         \filldraw[black] (\x,\y) circle[radius=0.5mm];}
        {\filldraw[black] (\x,\y) circle[radius=0.1mm];}
    }
  }
  \filldraw[black] (0,0) circle[radius=0.5mm];
\end{tikzpicture}
\end{center}
\caption{The visible points in $\Icc0\X\times\Icc0\Y$
  with the line connecting them to the origin}
\label{fig:visiblepoints_small}
\end{figure}%

The indicator function $\ifu\cV\from\Z^d\to\cA:=\{0,1\}$
can serve as a potential for a Schr\"odinger operator:
\begin{equation*}
  \Hvis\from\ell^2(\Z^d)\to\ell^2(\Z^d)\textq,
  \Hvis:=-\Laplace+\ifu\cV\text.
\end{equation*}
Analogous to \cref{sec:Anderson_Zd},
we can define the eigenvalue counting function
\begin{equation*}
  n\from\cF\to\B\textq,
  n(\Lambda)(E):=\tr\ifu{\Ioc{-\infty}E}(\Hvis^\Lambda)\text.
\end{equation*}
The results in \cite{LenzMV-08} imply that the thermodynamic limit
\begin{equation*}
  \lim_{L\to\infty}\setsize{\Lambda_L}^{-1}n(\Lambda_L)(E)
\end{equation*}
exists uniformly for all $E\in\R$,
but, to the best of our knowledge, there is no Pastur--Shubin formula.

\begin{figure}
\pgfmathsetmacro{\X}{150}
\pgfmathsetmacro{\Y}{100}
\begin{center}
\pgfmathsetmacro{\scaleforpicture}{12/(\X+1)}
\tikzsetnextfilename{visible2}%
\begin{tikzpicture}[scale=\scaleforpicture]
  \filldraw (0,0) circle[radius=0.5mm];
  \foreach \y in {0, ..., \Y} {
    \foreach \x in {0, ..., \X} {
      \pgfmathsetmacro{\visible}{ gcd(\x,\y) == 1 }
      \ifthenelse{\visible = 1}
        {\filldraw (\x,\y) circle[radius=0.5mm];}
        {}
    }
  }
\end{tikzpicture}
\end{center}
\caption{The visible points in $\Icc0\X\times\Icc0\Y$}
\label{fig:visiblepoints_large}
\end{figure}
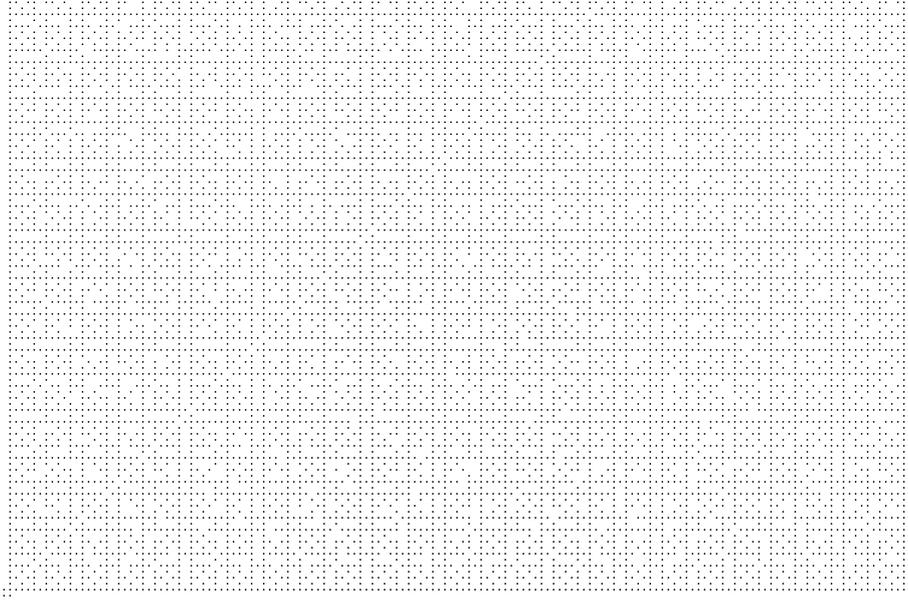%

\subsection{Patterns and frequencies}\label{sec:patternsfrequencies}

To understand the mechanism behind \cref{thm:LMV-ecf} better
and to motivate the abstract formulation of the next theorem,
let us examine the situation of percolation without potential in more detail.
We already indicated in \cref{sec:siteedge}
with the example of the eigenvalue~$0$, how discontinuities of~$N$ arise.
We now want to understand, how the limit of the eigenvalue counting functions
on large boxes obtains discontinuities of the same height as the IDS.
Let us focus on finite clusters again.
In $n(\Lambda_L,\omega)$,
the eigenvalues corresponding to eigenfunctions supported on finite clusters
are counted at least as often as a copy of their cluster
is contained in~$\Lambda_L$.
Because we normalize with the volume of the box,
$\setsize{\Lambda_L}^{-1}n(\Lambda_L,\omega)$,
the important quantity turns out to be the relative frequency
with which a finite cluster occurs in a large box~$\Lambda_L$.
Since the translations $(\tau_g)_{g\in G}$ act ergodically on~$\Omega$,
the relative frequencies converge to the probability
of their respective cluster at a fixed location.
But exactly these probabilities cause the discontinuities of the IDS.

Let us formalize the counting of copies of clusters, or,
in the more general setting of a finite set of colors~$\setsize\cA<\infty$,
shifts of patterns in boxes, following \cite{LenzMV-08}.
Recall that the set of finite subsets of~$G$ is denoted by~$\cF$.
A (finite) \emph{pattern} with \emph{domain} $\Lambda\in\cF$
is a map $P\from\Lambda\to\cA$.
Since a pattern $P\in\cA^\Lambda$ can be identified with the set
$\{\omega\in\Omega\mid\omega_\Lambda=P\}$,
we reuse the notation for the group action of~$G$ on~$\Omega$
for the shifts of patterns:
\begin{equation*}
  \tau_g\from\cA^\Lambda\to\cA^{\tilde\tau_g\Lambda}\textq,
  (\tau_gP)(v):=P(\tilde\tau_gv)\text.
\end{equation*}
This $G$-action induces an equivalence relation on the set of all finite patterns.
We denote the equivalence class of a pattern $P\from\Lambda\to\cA$
by $[P]_G:=\{\tau_gP\mid g\in G\}$.

A coloring $\omega\in\Omega=\cA^G$ and a finite subset $\Lambda\in\cF$
define the pattern $\omega_\Lambda:=\omega|_\Lambda\from\Lambda\to\cA$,
$\omega_\Lambda(v):=\omega_v$.
Similarly, for $\Lambda\subseteq\Lambda'\in\cF$, a pattern~$P'$ with domain~$\Lambda'$
induces the pattern $P'_\Lambda:=P'|_\Lambda\from\Lambda\to\cA$ via $P'_\Lambda(v):=P'(v)$.
Given $\Lambda'\in\cF$, set
\begin{equation*}
  P'|_\cF:=
  \{P'|_\Lambda\mid\Lambda\in\cF,\Lambda\subseteq\Lambda'\}\text.
\end{equation*}
This set of all finite patterns induced by $P'\in\cA^{\Lambda'}$
is useful to count how often a copy of a pattern $P\in\cA^\Lambda$
occurs in~$P'$:
\begin{equation}\label{eq:patterncount}
  \sharp_PP':=\setsize{[P]_G\isect P'|_\cF}\text.
\end{equation}
See Figure \ref{fig:patterncount} for a visualisation of this pattern counting function.
\begin{figure}
 
\begin{tikzpicture}[scale=0.6]
{
{\color{black}

{

\foreach \y in {0,1,...,8}{
\draw (1,\y) --(9,\y);
}
\foreach \x in {1,2,...,9}{
\draw (\x,0) --(\x,8);
}

\filldraw[blue] (8,2) circle (4pt);
\filldraw[blue] (7,8) circle (4pt);
\filldraw[blue] (6,5) circle (4pt);
\filldraw[blue] (7,2) circle (4pt);
\filldraw[blue] (6,8) circle (4pt);

\filldraw[red] (1,1) circle (4pt);
\filldraw[green] (1,2) circle (4pt);
\filldraw[green] (1,3) circle (4pt);
\filldraw[red] (1,4) circle (4pt);
\filldraw[red] (1,5) circle (4pt);
\filldraw[red] (1,6) circle (4pt);
\filldraw[green] (1,7) circle (4pt);
\filldraw[green] (1,8) circle (4pt);
\filldraw[red] (2,1) circle (4pt);
\filldraw[green] (2,1) circle (4pt);
\filldraw[red] (2,2) circle (4pt);
\filldraw[green] (2,3) circle (4pt);
\filldraw[blue] (2,4) circle (4pt);
\filldraw[blue] (2,5) circle (4pt);
\filldraw[red] (2,6) circle (4pt);
\filldraw[red] (2,7) circle (4pt);
\filldraw[green] (2,8) circle (4pt);
\filldraw[green] (3,1) circle (4pt);
\filldraw[red] (3,2) circle (4pt);
\filldraw[red] (3,3) circle (4pt);
\filldraw[red] (3,4) circle (4pt);
\filldraw[green] (3,5) circle (4pt);
\filldraw[red] (3,6) circle (4pt);
\filldraw[blue] (3,7) circle (4pt);
\filldraw[green] (3,8) circle (4pt);
\filldraw[blue] (4,1) circle (4pt);
\filldraw[green] (4,2) circle (4pt);
\filldraw[green] (4,3) circle (4pt);
\filldraw[green] (4,4) circle (4pt);
\filldraw[blue] (4,5) circle (4pt);
\filldraw[blue] (4,6) circle (4pt);
\filldraw[blue] (4,7) circle (4pt);
\filldraw[blue] (4,8) circle (4pt);
\filldraw[red] (5,1) circle (4pt);
\filldraw[blue] (5,2) circle (4pt);
\filldraw[red] (5,3) circle (4pt);
\filldraw[green] (5,4) circle (4pt);
\filldraw[green] (5,5) circle (4pt);
\filldraw[blue] (5,6) circle (4pt);
\filldraw[blue] (5,7) circle (4pt);
\filldraw[blue] (5,8) circle (4pt);
\filldraw[red] (6,1) circle (4pt);
\filldraw[red] (6,2) circle (4pt);
\filldraw[green] (6,3) circle (4pt);
\filldraw[red] (6,4) circle (4pt);
\filldraw[green] (6,5) circle (4pt);
\filldraw[red] (6,6) circle (4pt);
\filldraw[red] (6,7) circle (4pt);
\filldraw[green] (6,8) circle (4pt);
\filldraw[green] (7,1) circle (4pt);
\filldraw[blue] (7,2) circle (4pt);
\filldraw[red] (7,3) circle (4pt);
\filldraw[red] (7,4) circle (4pt);
\filldraw[green] (7,5) circle (4pt);
\filldraw[red] (7,6) circle (4pt);
\filldraw[green] (7,7) circle (4pt);
\filldraw[green] (7,8) circle (4pt);
\filldraw[green] (8,1) circle (4pt);
\filldraw[blue] (8,2) circle (4pt);
\filldraw[green] (8,3) circle (4pt);
\filldraw[green] (8,4) circle (4pt);
\filldraw[blue] (8,5) circle (4pt);
\filldraw[blue] (8,6) circle (4pt);
\filldraw[red] (8,7) circle (4pt);
\filldraw[green] (8,8) circle (4pt);
\filldraw[red] (9,1) circle (4pt);
\filldraw[blue] (9,2) circle (4pt);
\filldraw[red] (9,3) circle (4pt);
\filldraw[green] (9,4) circle (4pt);
\filldraw[blue] (9,5) circle (4pt);
\filldraw[blue] (9,6) circle (4pt);
\filldraw[blue] (9,7) circle (4pt);
\filldraw[red] (9,8) circle (4pt);
\filldraw[blue] (1,0) circle (4pt);
\filldraw[blue] (2,0) circle (4pt);
\filldraw[red] (3,0) circle (4pt);
\filldraw[red] (4,0) circle (4pt);
\filldraw[green] (5,0) circle (4pt);
\filldraw[blue] (6,0) circle (4pt);
\filldraw[green] (7,0) circle (4pt);
\filldraw[red] (8,0) circle (4pt);
\filldraw[green] (9,0) circle (4pt);

}
}
}

\draw[rounded corners=0.1cm]
 (4.7,4.2)-- ++(0,1.1)-- ++(1.6,0)-- ++(0,-1.6)-- ++(-1.6,0)-- ++(0,0.5);
\draw[rounded corners=0.1cm]
 (0.7,2.2)-- ++(0,1.1)-- ++(1.6,0)-- ++(0,-1.6)-- ++(-1.6,0)-- ++(0,0.5);
\draw[rounded corners=0.1cm]
 (6.7,0.2)-- ++(0,1.1)-- ++(1.6,0)-- ++(0,-1.6)-- ++(-1.6,0)-- ++(0,0.5);
\draw[rounded corners=0.1cm]
 (0.7,7.2)-- ++(0,1.1)-- ++(1.6,0)-- ++(0,-1.6)-- ++(-1.6,0)-- ++(0,0.5);
\draw[rounded corners=0.1cm]
 (3.7,3.2)-- ++(0,1.1)-- ++(1.6,0)-- ++(0,-1.6)-- ++(-1.6,0)-- ++(0,0.5);
\draw[rounded corners=0.1cm]
 (7.7,3.2)-- ++(0,1.1)-- ++(1.6,0)-- ++(0,-1.6)-- ++(-1.6,0)-- ++(0,0.5);
\draw[rounded corners=0.1cm]
 (6.7,7.2)-- ++(0,1.1)-- ++(1.6,0)-- ++(0,-1.6)-- ++(-1.6,0)-- ++(0,0.5);

\end{tikzpicture}
\caption{The figure shows a pattern $P'$ defined on a square $Q$ of side length~$9$. The marked $2\times2$ squares highlight the $7$ copies of a pattern $P$ in $P'$. Thus, we count $\sharp_PP'=7$. Obviously, to obtian a meaningful proportion of the occurrences of $P$ in $P'$, an appropriate normalization term is the size of the domain of $P'$, i.\,e.~$\tfrac{\sharp_PP'}{\setsize Q}=\tfrac{7}{81}$.}
\label{fig:patterncount}
\end{figure}
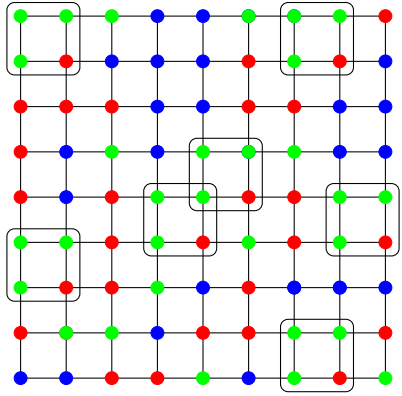

In this notation and for the probabilistic models from
\cref{sec:defsandmodels}, the ergodic theorem for $\Z^d$-actions,
see \cite{Keller-98}, states that, for all $\Lambda\in\cF$, $P\in\cA^\Lambda$,
and $\Prob$-almost all $\omega\in\Omega$,
\begin{equation*}
  \setsize{\Lambda_L}^{-1}\sharp_P(\omega_{\Lambda_L})
  \xto{L\to\infty}\Prob\{\omega\in\Omega\mid\omega_\Lambda=P\}\text.
\end{equation*}
The existence of the limit
$\lim_{L\to\infty}\setsize{\Lambda_L}^{-1}\sharp_P(\cV|_{\Lambda_L})$ 
for all patterns $P\in\cA^\Lambda$ with $\Lambda\in\cF$
in the setting of \cref{sec:visible}
is shown in \cite{BaakeMP-00} with different methods.
See \cref{fig:visiblepoints_large} for an optical impression.

Another feature to extract from the percolation example is the following.
By restricting the Hamiltonian to finite boxes, we modify some clusters
so that some part of their boundary is more ``straight''.
As a consequence, the clusters with one ``straight'' boundary
are over-represented in the sample.
Of course, as the boxes grow larger and larger,
their surface increases as well.
Fortunately, the even faster growth of the volume of the boxes makes
the boundary negligible, or more precisely:
the proportion of the surface to the volume vanishes in the limit.
By this mechanism, the surplus of clusters with artificial straight boundary
becomes negligible.

To formalize the notions above, let us introduce the \emph{$r$-boundary}
of~$\Lambda\in\cF$ for $r>0$ as
\begin{equation}\label{eq:rboundary}
  \partial^r\Lambda:=\{x\in G\setminus\Lambda\mid\dist(x,\Lambda)\le r\}
    \union\{v\in\Lambda\mid\dist(v,G\setminus\Lambda)\le r\}\text.
\end{equation}
The distance~$\dist$ is the length of the shortest path in the Cayley graph,
or, for $G=\Z^d$, the distance induced by the $1$-norm.
A sequence $(Q_j)_{j\in\N}$ of finite sets $Q_j\in\cF$
is called a \emph{F\o{}lner sequence}, if, for all $r>0$,
\begin{equation}\label{eq:folner}
  \lim_{j\to\infty}\frac{\setsize{\partial^rQ_j}}{\setsize{Q_j}}
  =0\text.
\end{equation}
Another common name for F\o{}lner sequences is \emph{van Hove sequence}.
The finite boxes $\Lambda_L:=\Ico0L^d\isect\Z^d\in\cF$, $L\in\N$,
form an example of a F\o{}lner sequence, because
$\setsize{\partial^r\Lambda_L}
  =(L+2\floor r)^d-(L-2\floor r)^d
  =4\floor r\sum_{k=0}^{d-1}(L-2\floor r)^k(L-2\floor r)^{d-1+k}
  \le4dr(L+2r)^{d-1}$
is of the order of~$L^{d-1}$, while $\setsize{\Lambda_L}=L^d$.
If, for a pattern $P\in\cA^\Lambda$, an $\omega\in\Omega$,
and a F\o{}lner sequence $(Q_j)_{j\in\N}$, the limit
\begin{equation*}
  \nu_P:=\lim_{j\to\infty}\frac{\sharp_P(\omega_{Q_j})}{\setsize{Q_j}}
\end{equation*}
exists, it is called the \emph{frequency of~$P$ along $(Q_j)_{j\in\N}$}.
See Figure \ref{fig:patterncount} for an illustration of one element of this sequence.
Note that the existence of $\nu_P$, for all patterns $P$ with finite domain, can be shown almost surely for all our examples
from \cref{sec:defsandmodels}. 

The eigenvalue counting functions are in a certain sense local enough
to reflect the effect of small boundary compared to the volume.
The following notions encapsulate the crucial properties.
\begin{Definition}
  A map $b\from\cF\to\Ico0\infty$ is called \emph{boundary term}, if
  \begin{itemize}[nosep]
  \item $b$ is invariant under~$G$: $b(\Lambda)=b(\tilde\tau_g\Lambda)$
    for all $g\in G$,
  \item $\lim_{j\to\infty}\setsize{Q_j}^{-1}b(Q_j)=0$
    for all F\o{}lner sequences $(Q_j)_{j\in\N}$,
  \item $b$ is bounded, i.\,e.\ its \emph{bound}
    $D_b:=\sup\{\setsize\Lambda^{-1}b(\Lambda)\mid
      \Lambda\in\cF,\Lambda\ne\emptyset\}<\infty$
    is finite, and
  \item for $\Lambda,\Lambda'\in\cF$ we have
    $b(\Lambda\union\Lambda')\le b(\Lambda)+b(\Lambda')$,
    $b(\Lambda\isect\Lambda')\le b(\Lambda)+b(\Lambda')$, and
    $b(\Lambda\setminus\Lambda')\le b(\Lambda)+b(\Lambda')$.
  \end{itemize}
\end{Definition}

\begin{figure}
\begin{center}
\tikzsetnextfilename{boundaryrules}%
\begin{tikzpicture}
  \draw [very thick] (0,0) -- (0.75,0.5)
    arc [start angle=30, end angle=335, x radius=2cm, y radius=1cm]
    -- cycle;
  \draw (0,0) -- (2,-1) -- (1.5,1) -- cycle;
  \draw (1,0) arc [ start angle=0
                  , end angle=360
                  , x radius=2cm
                  , y radius=1cm
                  ] -- cycle;
  \node at (-1.3,0) {$\Lambda\setminus\Lambda'$};
  \node at (-2.9,.8) {$\Lambda$};
  \node at (2,.6) {$\Lambda'$};
\end{tikzpicture}
\end{center}
\caption{The requirement
  $b(\Lambda\setminus\Lambda')\le b(\Lambda)+b(\Lambda')$
  is motivated by the fact that the boundary of $\Lambda\setminus\Lambda'$
  is a subset of the boundary of~$\Lambda$
  united with the boundary of~$\Lambda'$.}
\label{fig:pacman}
\end{figure}
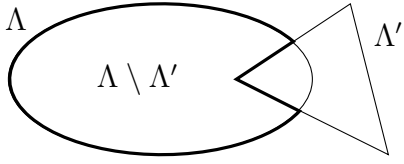
The last property is natural to require, as illustrated in \cref{fig:pacman},
but in fact only necessary when used with quasi tilings, see \cref{sec:general}.

\begin{Definition}\label{def:almost-additive}
  Consider a Banach space $(\B,\norm\argmt)$,
  a field $F\from\cF\to\B$, and a coloring $\omega\in\cA^G$.
  \begin{itemize}
  \item The field~$F$ is \emph{almost additive},
    if there exists a boundary term~$b$ such that,
    for all $n\in\N$, disjoint sets $\Lambda_1,\dotsc\Lambda_n\in\cF$,
    and $\Lambda:=\Union_{k=1}^n\Lambda_k$, it holds true that
    \begin{equation*}
      \norm{F(\Lambda)-\sum_{k=1}^nF(\Lambda_k)}
      \le\sum_{k=1}^nb(\Lambda_k)\text.
    \end{equation*}
  \item The field~$F$ is \emph{$\omega$-invariant},
    if for all $\Lambda,\Lambda'\in\cF$
    such that the patterns $\omega_\Lambda$ and $\omega_{\Lambda'}$
    are $G$-equivalent, we have
    \begin{equation*}
      F(\Lambda)=F(\Lambda')\text.
    \end{equation*}
  \end{itemize}
\end{Definition}

To an $\omega$-invariant field $F\from\cF\to\B$, the \emph{pattern function}
$\tilde F\from\Union_{\Lambda\in\cF}\cA^\Lambda\to\B$ with
\begin{equation*}
  \tilde F(P):=
  \begin{cases}
    F(\Lambda)&\text{if there is a set $\Lambda\in\cF$
      with $P\in[\omega_\Lambda]_G$, and}\\
    0&\text{otherwise}
  \end{cases}
\end{equation*}
is well defined.
Every almost additive field and $\omega$-invariant~$f$
is \emph{bounded} and has a \emph{bound} $C_F$ in the following sense:
\begin{equation*}
  C_F:=
  \sup\{\setsize\Lambda^{-1}\norm{F(\Lambda)}\mid
    \Lambda\in\cF\setminus\{\emptyset\}\}<\infty\text.
\end{equation*}
Indeed, since~$\cA$ is finite,
\begin{equation*}
  C_F
    \le\sup_{\Lambda\ne\emptyset}\frac1{\setsize\Lambda}\sum_{v\in\Lambda}
      \bigl(\norm{F(\{v\})}+b(\{v\})\bigr)
    \le\max_{a\in\cA^{\{0\}}}\tilde F(a)+b(\{0\})
    <\infty
\end{equation*}

The eigenvalue counting functions are good examples for these notions.
\begin{Proposition}[{\cite[Proposition~2]{LenzMV-08}}]
  For almost all $\omega\in\Omega$ and
  with respect to the Banach space $(\B,\norm[\infty]\argmt)$
  of right continuous bounded $\R$-valued functions on~$\R$,
  the eigenvalue counting functions~$n$ of the models given in
  \cref{sec:defsandmodels} 
  are $\omega$-invariant, almost additive with boundary term
  $b(\Lambda):=4\setsize{\partial^1\Lambda}$, $D_b\le8d+4$,
  and bounded with bound~$C_{n(\argmt,\omega)}=1$.
\end{Proposition}

The following theorem thus applies
to the eigenvalue counting functions of \cref{sec:defsandmodels}.
We emphasize that the error estimates imply that the IDS is approximated
by the eigenvalue counting functions uniformly in the energy.

\begin{Theorem}[{\cite[Theorem~1]{LenzMV-08}}]\label{thm:LenzMV}
  Let $\cA$ be a finite set, $\omega\in\cA^G$ a coloring,
  $(\B,\norm\argmt)$ an arbitrary Banach space,
  $(Q_j)_{j\in\N}$ a F\o{}lner sequence,
  $F\from\cF\to\B$ a bounded, $\omega$-invariant, and almost additive field
  with bound~$C_F$, pattern function $\tilde F$, boundary term~$b$,
  and bound~$D_b$ of~$b$.
  \par
  Assume that for every finite pattern~$P\from\Lambda\to\cA$, $\Lambda\in\cF$,
  the frequency $\nu_P:=\lim_{j\to\infty}\setsize{Q_j}^{-1}
    \sharp_P(\omega|_{Q_j})$ exists.
  Then, the limits
  \begin{equation}\label{eq:Fbar-F-Ftilde}
    \Fbar
    :=\lim_{j\to\infty}\frac{F(Q_j)}{\setsize{Q_j}}
    =\lim_{L\to\infty}\sum_{P\in\cA^{\Lambda_L}}
      \nu_P\frac{\tilde F(P)}{\setsize{\Lambda_L}}
  \end{equation}
  exist and are equal.
  Moreover, for all $j,L\in\N$, the bounds
  \begin{equation}\label{eq:Fbar-Ftilde}
    \biggnorm{\Fbar-\sum_{P\in\cA^{\Lambda_L}}
      \nu_P\frac{\tilde F(P)}{\setsize{\Lambda_L}}}
    \le\frac{b(\Lambda_L)}{\setsize{\Lambda_L}}
  \end{equation}
  and
  \begin{multline}\label{eq:Ftilde-F}
    \biggnorm{\frac{F(Q_j)}{\setsize{Q_j}}
      -\sum_{P\in\cA^{\Lambda_L}}
      \nu_P\frac{\tilde F(P)}{\setsize{\Lambda_L}}}\\
    \le\frac{b(\Lambda_L)}{\setsize{\Lambda_L}}
      +(C_F+D_b)\frac{\setsize{\partial^LQ_j}}{\setsize{Q_j}}
      +C_F\sum_{P\in\cA^{\Lambda_L}}
        \Bigabs{\frac{\sharp_P(\omega|_{Q_j})}{\setsize{Q_j}}-\nu_P}
  \end{multline}
  hold true.
\end{Theorem}

The error estimates are the crucial part of the theorem.
We note that there are two length scales, indexed by~$j$ and~$L$,
in the approximation of the limiting object~$\Fbar$.
They correspond to two stages of approximation in the proof.
In the first step,~$\Fbar$ is compared to the weighted average of
the contributions of the patterns on~$\Lambda_L$,
the weights being the frequencies of the patterns.
According to~\eqref{eq:Fbar-Ftilde}, the patterns on~$\Lambda_L$
capture the behavior of the limiting function up to an error of
size $\setsize{\Lambda_L}^{-1}b(\Lambda_L)$.
We mentioned above, that $(\Lambda_L)_L$ is a F\o{}lner sequence,
so by choosing~$L$ large enough, we can make this error term small.

\Cref{eq:Ftilde-F} addresses the problem that we used the limiting frequencies
in~\eqref{eq:Fbar-F-Ftilde} and~\eqref{eq:Fbar-Ftilde}
instead of the actual number of occurrences of patterns
in the finite region~$Q_j$.
This error bound requires us to choose~$j$
so large that the empirical frequencies of all patterns on~$\Lambda_L$
get close to their actual frequencies.
To recapitulate: First, we have to choose~$L$
large enough to make the patterns on~$\Lambda_L$
meaningful for the limit.
Then we have to choose~$j$ large enough such that the patterns on~$\Lambda_L$
are actually observed in proportions
that are close to the asymptotic frequencies.

In the random Schr\"odinger operator settings from \cref{sec:defsandmodels},
the second error is controlled by the randomness.
In fact, the ergodic theorem predicts that
the relative number of occurrences of a pattern
converges to its frequencies almost surely.
It follows from the theory of large deviations
that the probability of a fixed difference between the two
is exponentially small for large samples, that is, for large~$j$.

The frequencies $(\nu_P)_{P\in\cA^{\Lambda_L}}$
as well as their empirical counterparts
$(\setsize{Q_j}^{-1}\sharp_P(\omega|_{Q_j}))_{P\in\cA^{\Lambda_L}}$
can be interpreted as a probability mass function on~$\cA^{\Lambda_L}$.
The $1$-norm of their difference in~\eqref{eq:Ftilde-F}
is also known as the \emph{total variation norm}\label{tv_mentioned}
of the corresponding probability measures.

In the remainder of this note, we present two generalizations
which correspond to the two types of errors described above.
The main property of the group~$G=\Z^d$ used in \cref{thm:LenzMV}
is that it hosts F\o{}lner sequences.
That this is actually the crucial property necessary for this proof
is made explicit by generalizing the result to amenable groups,
which are characterized by the existence of a F\o{}lner sequence,
see \cref{sec:amenable}.

The second generalization concerns the finiteness of the set of colors~$\cA$.
The main obstacle to overcome this restriction is the probabilistic error.
The total variation norm of the difference of the distribution
and the empirical measure of a random variable
does not converge to zero for continuous random variables.
In order to allow infinitely many colors,
we will be forced to exploit more properties
of the eigenvalue counting functions,
see \cref{sec:infiniteA}.

\section{Ergodic theorems for finite colors on amenable groups}\label{sec:amenable}

In this section we discuss how the above ideas generalize to less restricted geometries.
In particular, we present Banach space-valued ergodic theorems for Cayley graphs generated by amenable groups.
Let us emphasize that the groups considered in this paper will always be finitely generated and therefore countable.
An amenable group is by definition a group containing subsets with an arbitrary small ratio between boundary and volume.
It is well known, that amenability is equivalent to the existence of a F\o{}lner sequence.

It turns out that even though amenability is the natural condition to generalize the geometry of $\Z^d$,
there is an additional requirement needed to almost directly implement the $\Z^d$-methods to this setting.
Here we are speaking about a so-called tiling condition, namely the condition that there exists a F\o{}lner sequence consisting of monotiles.
Obviously, in~$\Z^d$ a sequence of cubes serves as such a sequence, since for each~$j$ the group~$\Z^d$ can be tiled with cubes of side length~$j$.
For an arbitrary amenable group it is not known whether such a sequence exists or not.
Therefore, this section is structured in a first part discussing the monotile situation and
a second more involved part dealing with the general amenable groups using the technique of quasi tilings.

We proceed with some definitions which generalize the notion of previous
sections. Given a finitely generated group $G$ with a finite and symmetric
generating set $S\subseteq G$, i.\,e.\ $G=\langle S\rangle$, $\setsize S<\infty$
and $S=S^{-1}\not\ni \mbox{id}$, the corresponding Cayley graph~$\Gamma$ has vertex set~$G$ and
two vertices $v,w\in G$ are connected if and only if $vs=w$ for some $s\in S$.
The induced graph distance, sometimes called word metric, is denoted by~$\dist$.
A sequence $(Q_j)$ of finite subsets of~$G$ is called a F\o{}lner sequence if
$\setsize{\partial^r Q_j}/\setsize{Q_j}\to 0$ as
$j\to\infty$ for all $r>0$, see \eqref{eq:folner}. Here the $r$-boundary
is given as in \eqref{eq:rboundary}.

For an amenable group~$G$, the introduction of the physical models of
\cref{sec:Anderson_Zd,sec:percolation} works completely
analogous. Due to the fact that for general groups the group action is
usually written as a multiplication, the only difference is that the
group action $\tilde\tau$ is defined here as
\begin{equation}\label{def:tau}
  \tilde\tau \from G\times G \to G ,\quad (g,v)\mapsto \tilde\tau_gv := vg^{-1}.
\end{equation}
As in \eqref{eq:schroedinger} the Sch\"odinger operator in the Anderson model is given by
\[
  H_\omega:=-\Laplace+V_\omega\from\ell^2(G)\to\ell^2(G)
\]
with Laplace operator $\Delta\from\ell^2 (G)\to\ell^2(G)$
as in~\eqref{eq:laplaceop} and a (random) potential $V\from\Omega\times G\to \cA$
as in~\eqref{eq:randompotential}. As before we consider the probability
space $(\Omega,\cB,\Prob) :=(\cA^G,\Borel(A)^{\tensor G},\Prob_0^{\tensor G})$
with a finite set $\cA$.
An element $\omega\in\Omega$ is then interpreted as a coloring of the elements
of the group using the finite set of colors~$\cA$.

Besides this, also the definition of site and edge percolation does not
depend on the $\Z^d$-structure and generalizes straightforwardly to the
case of amenable groups. Thereby the Anderson model on percolation Cayley
graphs over finitely generated amenable groups is well-defined.

\subsection{Symmetrical tiling condition}
As mentioned before an additional condition is needed in order to implement the methods of~$\Z^d$ to amenable groups.
This so-called symmetric tiling condition is formulated as follows.
\begin{Definition}\label{def:STamenable}
 Let $G$ be a group. A subset $\Lambda\subseteq G$ \emph{symmetrically tiles} $G$ if there exists a set $T$ such that
 \begin{enumerate}[(i)]
  \item $T=T^{-1}$,
  \item $\Lambda T = G$,
  \item $\{\Lambda t \mid t\in T \}$ are pairwise disjoint.
 \end{enumerate}
 An amenable $G$ is said to satisfy the \emph{symmetric tiling condition},
 if there exists a F\o{}lner sequence $(\Lambda_L)$ such that each
 $\Lambda_L$, $L\in \N$, symmetrically tiles $G$. In this situation we
 call $G$ an \emph{ST-amenable group}.
\end{Definition}
Note that here $\Lambda T:=\{xt \mid x\in \Lambda, t\in T\}$ and $\Lambda t=\{xt\mid x\in \Lambda\}$.

Let us briefly remark on this condition.
Krieger proved in \cite{Krieger-07} based on work of Weiss \cite{Weiss-01}
that an amenable group satisfies the symmetrical tiling condition if it
is residually finite.
For instance, each group of polynomial volume growth is nilpotent
(by Gromov's theorem) and thus residually finite.
Since it is also of subexponential growth, it is amenable, too,
and hence ST-amenable.

An intensively studied and slightly more general condition than the
one stated in \cref{def:STamenable} can be obtained when not
assuming the symmetry (i). In this situation a set $Q$ satisfying (ii)
and (iii) is usually referred to as a \emph{monotile} and a group $G$
containing a F\o{}lner sequence consisting only of monotiles is called
\emph{monotileable}.
Let us remark that it is still not known if there exists an amenable
group which is not monotileable.

In the situation of ST-amenable groups one can prove the following:
\begin{Theorem}[{\cite{LenzSV-10,LenzSV-12}}]\label{thm:LenzSV}
 Let $G$ be a finitely generated ST-amenable group. Let~$(Q_j)$
 and~$(\Lambda_L)$ be F\o{}lner sequences such that each $\Lambda_L$,
 $L\in\N$ symmetrically tiles~$G$.
 Let $\cA$ be finite and $\omega\in \cA^G$ be given.

 Moreover, let $(\B,\norm\argmt)$ be a Banach space and $F\from\cF\to\B$
 a bounded, $\omega$-invariant and almost additive field with bound~$C_F$,
 pattern function $\tilde F$, boundary term~$b$,
  and bound~$D_b$ of~$b$. As before we use the notation $\cF=\{A\subseteq G\mid A \text{ finite}\}$.

 Assume that for every finite pattern~$P\from\Lambda\to\cA$, $\Lambda\in\cF$,
  the frequency $\nu_P:=\lim_{j\to\infty}\setsize{Q_j}^{-1}
    \sharp_P(\omega|_{Q_j})$ exists.
  Then, the limits
  \begin{equation}\label{eq:Fbar-F-Ftilde.amenable}
    \Fbar
    :=\lim_{j\to\infty}\frac{F(Q_j)}{\setsize{Q_j}}
    =\lim_{L\to\infty}\sum_{P\in\cA^{\Lambda_L}}
      \nu_P\frac{\tilde F(P)}{\setsize{\Lambda_L}}
  \end{equation}
  exist and are equal.
  Moreover, for all $j,L\in\N$, the bound
  \begin{multline}\label{eq:Ftilde-F.amenable}
    \biggnorm{\frac{F(Q_j)}{\setsize{Q_j}}
      -\sum_{P\in\cA^{\Lambda_L}}
      \nu_P\frac{\tilde F(P)}{\setsize{\Lambda_L}}}\\
    \le\frac{b(\Lambda_L)}{\setsize{\Lambda_L}}
      +(C_F+D_b)\frac{\setsize{\partial^{\diam(\Lambda_L)}Q_j}}{\setsize{Q_j}}
      +C_F\sum_{P\in\cA^{\Lambda_L}}
        \Bigabs{\frac{\sharp_P(\omega|_{Q_j})}{\setsize{Q_j}}-\nu_P}
  \end{multline}
  holds true. Here, $\diam(\Lambda_L)=\max\{\dist(x,y)\mid x,y\in \Lambda_L\}$ is the diameter of the finite set~$\Lambda_L$.
\end{Theorem}

When comparing \cref{thm:LenzMV} and \cref{thm:LenzSV},
it turns out that the transition from $\Z^d$ to ST-amenable groups does not imply a quantitative difference in the strength of the result.
In particular, the only difference is that in the $\Z^d$ setting $(\Lambda_L)$ is a sequence of cubes with side length $L$ (which naturally tiles $\Z^d$) and in the ST-amenable group setting  $(\Lambda_L)$ is a F\o{}lner sequence assumed to be symmetrically tiling.
In the error bound this results in the substitution of the side length $L$ of a cube by the diameter of $\Lambda_L$.

In the following we give a sufficient condition for the existence of the frequencies in the situation of a random coloring.
Here we need the notion of a \emph{tempered} F\o{}lner sequence, which is a F\o{}lner $(Q_j)$ sequence with the additional property that there is some $C>0$ such that
\[
  \Biggl|\bigcup_{k<n}Q_k^{-1}Q_j\Biggr| \leq C \setsize{Q_j}
\]
for all $j\in\N$.
Each F\o{}lner sequence has a tempered subsequence, see \cite{Lindenstrauss-01}.

\begin{Theorem}
  Let~$G$ be a finitely generated amenable group,
  $\cA$ some finite set and let~$\mu$ be a probability measure on
  $(\Omega,\cB) =(\cA^G,\Borel(A)^{\tensor G})$.
  We assume that the action $\tilde\tau$ given via~\eqref{def:tau}
  of~$G$ on~$\Omega$ is measure preserving and ergodic w.\,r.\,t.~$\mu$.
  Then, for any tempered F\o{}lner sequence $(Q_j)$,
  there exists a event $\tilde\Omega$ of full measure, such that the limit
  \begin{equation*}
    \lim_{j\to\infty}\frac{\sharp_P(\omega|_{Q_j})}{\abs{Q_j}}
  \end{equation*}
  exists for all patterns $P\in\bigcup_{Q\in\cF}\cA^{Q}$
  and all $\omega\in\tilde\Omega$.
  Moreover, the limit is deterministic in the sense that it is
  independent of the specific choice of~$\omega\in\tilde\Omega$.
\end{Theorem}


The above result is a direct consequence (see for instance \cite{Schwarzenberger-08})
of Lindenstrauss ergodic theorem \cite{Lindenstrauss-01}.
We omit the precise definitions of measure preserving and ergodic action.
However, we want to emphasize that in the particular case where
$\mu=\Prob_0^{\tensor G}$ is a product measure,
the assumptions on~$\tilde\tau$ are met.
Thus, in the percolation setting of \cref{sec:percolation}
one obtains that almost all configurations satisfy the assumption
of well-defined frequencies.


\subsection{General amenable groups}\label{sec:general}

As outlined before, it is still not known if each amenable group satisfies the symmetrical tiling condition.
Roughly speaking, in the previous sections this tiling condition is the crucial tool to mediate between the two F\o{}lner sequence $(\Lambda_L)$ and $(Q_j)$.
Thus, when considering amenable groups without an additional tiling assumption the situation is far more challenging.

A way to overcome this lack is to apply the theory of $\epsilon$-quasi tilings developed by Ornstein and Weiss \cite{OrnsteinW-87} in 1987, see also \cite{PogorzelskiS-16} for the quantitative estimates used is the present setting.
The key idea here is to soften the condition of a \emph{perfect} tiling with copies of \emph{one} set taken from a F\o{}lner sequence, and rather
\begin{itemize}
  \item use finitely many different sets of the F\o{}lner sequence, and
  \item allow imperfectness in the tiling (in sense of small \emph{overlaps} and \emph{uncovered areas}).
\end{itemize}

More precisely we use the following definition:
\begin{Definition}\label{def:quasitiling}
 Let $G$ be a finitely generated group, $Q\subseteq G$ a finite set and $\epsilon>0$.
 We say that $K_1,\dots,K_N\subseteq G$ with center sets $T_1,\dots, T_N$, short $(K_i,T_i)_{i=1}^N$, are an \emph{$\epsilon$-quasi tiling} of~$Q$ if
 \begin{enumerate}[(i)]
  \item the sets $K_iT_i$, $i\in\{1,\dotsc,N\}$ are pairwise disjoint and subsets of $Q$;
  \item $\setsize{Q\setminus \bigcup_{i=1}^N K_iT_i} \leq 2\epsilon \setsize{Q}$;
  \item there are subsets $\mathring{K_i}\subseteq K_i$, $i\in\{1,\dotsc,N\}$, such that
  \begin{itemize}
    \item for each $i$ the sets $\mathring{K_i} t$, $t\in T_i$ are pairwise disjoint
    \item $\setsize{K_i\setminus \mathring{K_i}}\leq \epsilon \setsize{K_i}$
  \end{itemize}
 \end{enumerate}
\end{Definition}

For technical reasons, the sequence $(\Lambda_L)$ which will provide the F\o{}lner sets $K_i$ to quasi tile the group is assumed to be \emph{nested},
i.\,e.\ for each $L\in\N$ we have $\id \in \Lambda_L\subseteq \Lambda_{L+1}$.
Note that, starting from an arbitrary F\o{}lner sequence,
one can construct a nested F\o{}lner sequence by translating elements of an appropriate subsequence.

In the following it turns  out to be convenient to define for given $\epsilon\in(0,1)$ and $i\in\N_0$ the numbers $N(\epsilon)$ and $\eta_i(\epsilon)$ by 
\begin{equation}\label{eq:def:N,eta}
	N(\epsilon):=\biggceil{\frac{\ln(\epsilon)}{\ln(1-\epsilon)}}
	\qtextq{and}
	\eta_i(\epsilon):=\epsilon(1-\epsilon)^{N(\epsilon)-i}
	\text.
\end{equation}
As usual we use the Gau\ss{}ian bracket notation $\ceil b:=\inf\{z\in\Z\mid z\ge b\}$.
The following theorem shows that~$N(\epsilon)$ is the number of required shapes~$K_i$ in order to $\epsilon$-quasi tile a set~$Q$.
Moreover, for fixed~$i$ the~$\eta_i(\epsilon)$ can be interpreted as the ratio of the points covered by copies of~$K_i$ in the $\epsilon$-quasi tiling.

\begin{Theorem}\label{thm:STP}
  Let~$G$ be a finitely generated amenable group,
	$(\Lambda_L)$ a nested F{\o}lner sequence, and $\epsilon\in\Ioo0{0.1}$.
  Then there is a finite and strictly increasing selection of sets
  $K_i\in\{\Lambda_L\mid L\in\N\}$, $i\in\{1,\dotsc,N(\epsilon)\}$,
	with the following property:
  For each F{\o}lner sequence~$(Q_j)$,
  there exists $j_0(\epsilon)\in\N$
  satisfying that for all $j\geq j_0(\epsilon)$ there exists sets~$T_1^j, \ldots, T_{N(\epsilon)}^j$ such that
  $(K_i,T_i^j)_{i=1}^{N(\epsilon)}$ is an $\epsilon$-quasi tiling of~$Q_j$.
	Moreover, for all $j\ge j_0(\ve)$ and all $i\in\{1,\dots,N(\epsilon)\}$,
	the proportion of~$Q_j$ covered by the tile~$K_i$ satisfies
	\begin{equation}\label{eq:etdensity}
		\biggabs{\frac{\setsize{K_iT_i^j}}{\setsize{Q_j}}-\eta_i(\ve)}
		\le\frac{\ve^2}{N(\ve)}\text.
	\end{equation}
\end{Theorem}

The proof of \cref{thm:STP} is to be found in \cite{PogorzelskiS-16}.
Although \cref{thm:STP} provides
all the elements used to formulate the desired ergodic theorem for
general amenable groups (\cref{thm:PSerg}), a far more involved
result is applied in the proof. More precisely, in the proof of
\cref{thm:PSerg} one does not only need \emph{one} possibility to quasi
tile a given set $Q_j$ with $K_1,\dots,K_{N(\epsilon)}$, but one rather
needs a bunch of \emph{different} possibilities to quasi tile $Q_j$
with $K_1,\dots,K_{N(\epsilon)}$. These different tilings of $Q_j$ need
to be chosen such that for (almost) all $g\in Q_j$ the frequency (over
different tilings) that it is covered by one $K_i$ is up to a small error
$\eta_i(\epsilon)/|K_i|$. In this sense, this covering result is referred
to as uniform $\epsilon$-quasi tiling.
We refer to \cite{PogorzelskiS-16} for details.

Let us formulate the ergodic theorem based on the $\epsilon$-quasi tiling results.

\begin{Theorem}\label{thm:PSerg} Assume:
\begin{itemize}
 \item $G$ is a finitely generated amenable group.
 \item $\cA$ is a finite set and $\omega\in \cA^G$.
 \item $(Q_j)$ is a F\o{}lner sequence such that
 the frequency $\nu_P:=\lim_{j\to\infty}\setsize{Q_j}^{-1}
    \sharp_P(\omega|_{Q_j})$ exists for every finite pattern~$P\from\Lambda\to\cA$, $\Lambda\in\cF$.
 \item $(\Lambda_L)$ is a nested F\o{}lner sequence.
 \item For given $\epsilon\in(0,\tfrac1{10})$ the sets $K_i$, $i\in\{1,\dots,N(\epsilon)\}$ are chosen according to \cref{thm:STP}.
 \item $(\B,\norm\argmt)$ is a Banach space and $F\from\cF\to\B$ is $\omega$-invariant, and almost additive with bound~$C_F$, pattern function $\tilde F$, boundary term~$b$,
 and bound~$D_b$ of~$b$.
\end{itemize}
 Then, the limits
 \begin{equation}
  \Fbar:=  \lim_{j\to\infty}\frac{F(Q_j)}{\abs{Q_j}} = \lim_{\substack{\epsilon \searrow 0 \\ \epsilon<0.1}} \sum_{i=1}^{N(\epsilon)} \eta_i(\epsilon) \sum_{P\in \cA^{K_i}} \nu_P \frac{\tilde F(P)}{\abs{K_i}}
 \end{equation}
 exist and are equal.
 Moreover, for given $\epsilon\in(0,\tfrac1{10})$ there exist some $j(\epsilon),r(\epsilon)\in\N$ such that for all $j\geq j(\epsilon)$ the bound
 \begin{multline}\label{eq:Ftilde-F.amenable.gen}
    \biggnorm{ \frac{F(Q_j)}{\abs{Q_j}} - \sum_{i=1}^{N(\epsilon)} \eta_i(\epsilon) \sum_{P\in \cA^{K_i}} \nu_P \frac{\tilde F(P)}{\abs{K_i}}}\\
    \le
   4 \sum_{i=1}^{N(\epsilon)}  \eta_i(\epsilon) \frac{b(K_i )}{|K_i|}      +(C_F+4D_b)\frac{\abs{\partial^{r(\epsilon)}Q_j}}{\abs{Q_j}}\sum_{i=1}^{N(\epsilon)} \abs{K_i}     \\
      +     C_F  \sum_{i=1}^{N(\epsilon)} \eta_i(\epsilon) \sum_{P \in \cA^{K_i}} \left| \frac{\sharp_P(\omega|_{Q_j})}{|Q_j|} - \nu_P \right|  + (11C_F+32D_b) \epsilon
  \end{multline}
  holds true.
\end{Theorem}

When comparing the estimate with the one in \cref{thm:LenzSV},
note that the difference \eqref{eq:Ftilde-F.amenable.gen} gets small if one firstly executes the
limit ${j\to\infty}$ and afterwards ${\epsilon\searrow 0}$. Here
the limit $\epsilon\searrow0$ corresponds to 
$L\to\infty$
in the
previous setting. For a detailed discussion why the error terms tend
to zero we refer to \cite{PogorzelskiS-16}. We confine ourselves to
the (rough) statement that the first three terms in the estimate
\eqref{eq:Ftilde-F.amenable.gen}  correspond in this ordering to
three terms in \eqref{eq:Ftilde-F.amenable} or \eqref{eq:Ftilde-F},
respectively.

\section{Glivenko--Cantelli type theorems}\label{sec:infiniteA}

In this section we will consider the situation that the set of colors~$\cA$ is no longer finite.
This means that we are leaving a combinatorial setting and relying on a probabilistic framework instead.
This has been already introduced for a number of examples  in \cref{sec:Anderson_Zd,sec:percolation}.
In particular, we will need some independence and monotonicity assumptions.
The monotonicity property will allow us to smooth out and regularize certain quantities which we otherwise do not know how to estimate.
The independence assumption gives us a tool to describe the existence of frequencies used in Section \ref{sec:finiteA}
in a constructive and quantitative manner.

We reformulate the notions of \cref{sec:finiteA}
involving fields $F\from\cF\to\B$ with values in an arbitrary Banach space~$\B$
for fields of the form $f\from\cF\times\Omega\to\B$
with a probability space~$\Omega$ and the specific Banach space
of right continuous functions with $\sup$-norm.
For easy distinction, we denote the latter with lowercase letters.
Almost additive is assumed to hold true uniformly on~$\Omega$.
The property $\omega$-invariance for fixed $\omega\in\Omega$
is split up in equivariance and locality.
\begin{Definition}\label{def:admissible}
  A field $f\from\cF\times\Omega\to\B$
  is
  \begin{itemize}
  \item
    \emph{almost additive}, if there is a boundary term
    $b\from\cF\to\Ico0\infty$
    such that for all $\omega\in\Omega$,
    pairwise disjoint $\Lambda_1,\dots,\Lambda_n\in\cF$,
    and $\Lambda:=\bigcup_{i=1}^n\Lambda_i$, we have
    \begin{equation*}
      \Bignorm{f(\Lambda,\omega)-\sum_{i=1}^nf(\Lambda_i,\omega)}
        \le\sum_{i=1}^nb(\Lambda_i)\text.
    \end{equation*}
  \item
    \emph{equivariant}, if for $\Lambda\in\cF$, $g\in G$ and $\omega\in\Omega$ we have
    $f(\tilde\tau_g\Lambda,\omega)=f(\Lambda,\tau_g\omega)$.
  \item
    \emph{local}, if for all $\Lambda\in\cF$ and $\omega,\omega'\in\Omega$,
    $\omega_\Lambda=\omega'_\Lambda$ implies
    $f(\Lambda,\omega)=f(\Lambda,\omega')$.
  \item
    \emph{bounded}, if
    $
      \sup_{\omega\in\Omega}\norm{f(\{\id\},\omega)}
      <\infty\text.
    $
  \end{itemize}
\end{Definition}


\subsection{Glivenko--Cantelli theory}

To simplify the motivation, we first consider $G=\Z^d$
and assume that the probability measure~$\Prob=\Prob_0^{\tensor G}$
on~$\Omega$ is a product measure.
This implies in particular that, for every local field
$f\from\cF\times\Omega\to\B$, the random variables
$\omega\mapsto f(\{0\},\tau_g\omega)$, $g\in G$, are independent.
Here, we used that for local~$f$, in particular,
$f(\{0\},\omega)$ depends only on~$\omega_0$
and not on~$\omega_{G\setminus\{0\}}$.

To explain the relation of our methods to Glivenko--Cantelli theory,
assume for the moment that $f\from\cF\times\Omega\to\B$
is not only almost additive but \emph{exactly} additive.
That means that for a finite subset $\Lambda\in\cF$ of~$G$
and a realization of colors $\omega\in\Omega=\cA^G$,
we can split $f(\Lambda,\omega)$
without any errors into a sum over singleton sets
\begin{equation*}
  f(\Lambda,\omega)
  =\sum_{v\in\Lambda}f(\{v\},\omega)\text.
\end{equation*}
By equivariance, we can rewrite $f(\{v\},\omega)=f(\{0\},\tau_v^{-1}\omega)$.
For the special case $\B=\R$,
the law of large numbers allows us to calculate the thermodynamic limit
\begin{equation*}
  \lim_{L\to\infty}\frac{f(\Lambda_L,\omega)}{\setsize{\Lambda_L}}
  =\lim_{L\to\infty}\frac1{\setsize{\Lambda_L}}
    \sum_{v\in\Lambda_L}f(\{0\},\tau_v^{-1}\omega)
  =\E[f(\{0\},\argmt)]
\end{equation*}
$\Prob$-almost surely.

Of course, in view of our examples in \cref{sec:defsandmodels},
we are more interested in the Banach space~$\B$
of right continuous functions from~$\R$ to~$\R$ with $\sup$-norm.
To head in this direction it is advantageous to lift the point of view to the empirical probability
$\widehat\Prob_\Lambda^\omega
  :=\setsize\Lambda^{-1}\sum_{v\in\Lambda}\delta_{\omega_v}$
and to write integration as dual pair $\spr\argmt\argmt$.
In this notation, the average from above can be written as
\begin{equation*}
  \setsize\Lambda^{-1}f(\Lambda,\omega)
  =\spr{f(\{0\},\argmt)}{\widehat\Prob_\Lambda^\omega}
  \text.
\end{equation*}
Here is a special case,
where a theorem from classical probability theory helps.
Assume that $\cA\in\Borel(\R)$, and let $f\from\cF\times\Omega\to\B$
count the number of random variables in~$\Lambda$
with value less than a given threshold $E\in\R$:
\begin{equation*}
  f(\Lambda,\omega)(E)
  :=\sum_{v\in\Lambda}\ifu{\Ico{\omega_v}\infty}(E)
  =\sum_{v\in\Lambda}\ifu{\Ioc{-\infty}E}(\omega_v)
  =\setsize\Lambda\spr{\ifu{\Ioc{-\infty}E}}{\widehat\Prob_\Lambda^\omega}\text.
\end{equation*}
This defines an additive, local, and equivariant field,
which can be interpreted as the (not normalized)
empirical distribution function of the sample~%
$(\omega_v)_{v\in\Lambda}$.
The thermodynamic limit in~$\B$, i.\,e.\ uniformly, is $\Prob$-almost surely
\begin{equation*}
  \lim_{L\to\infty}\setsize{\Lambda_L}^{-1}f(\Lambda_L,\omega)
  =\lim_{L\to\infty}
    \spr{\ifu{\Ioc{-\infty}\argmt}}{\widehat\Prob_\Lambda^\omega}
  =\spr{\ifu{\Ioc{-\infty}\argmt}}{\Prob_0}
  =\Prob(\omega_0\le\argmt)
\end{equation*}
by the theorem of Glivenko and Cantelli:

\begin{Theorem}[\cite{Glivenko-33,Cantelli-33}]\label{thm:GC}
  Let $V_j$, $j\in\N$, be real valued, independent
  and identically distributed random variables on
  $(\Omega,\Borel,\Prob)$ and $\widehat\Prob_n:=\frac1n\sum_{j=1}^n\delta_{V_j}$
  the corresponding empirical distribution.
  Then there exists an event $\Omega_{\mathrm{unif}}$
  with probability $\Prob(\Omega_{\mathrm{unif}})=1$
  such that for all $\omega\in\Omega_{\mathrm{unif}}$:
  \begin{equation*}
    \sup_{E\in\R}\abs{\spr{\ifu{\Ioc{-\infty}E}}{\widehat\Prob_n(\omega)-\Prob}}
    \xto{n\to\infty}0\text.
  \end{equation*}
\end{Theorem}

We learn that the choice of this Banach space~$\B$
means to prove the convergence of~$\widehat\Prob_\Lambda^\omega$ to~$\Prob_0$
with respect to a supremum over appropriate test functions.
The route pursued in \cref{thm:LenzMV,thm:LenzSV,thm:PSerg}
for finite alphabets corresponds to the estimate
\begin{equation*}
  \abs{\spr g{\widehat\Prob_\Lambda^\omega-\Prob_0}}
  \le\norm[\infty]g\TVnorm{\widehat\Prob_\Lambda^\omega-\Prob_0}
\end{equation*}
with $g:=f(\{0\},\argmt)$.
But, as the next example shows, for smooth random variables,
the difference does not converge to zero in total variation.

\begin{Example}
  Let $\cA=\Icc01$ and~$X_n$, $n\in\N$,
  be real valued i.\,i.\,d.\ random variables,
  distributed uniformly on~$\cA$.
  The empirical distribution
  $\widehat\Prob_n:=n^{-1}\sum_{j=1}^n\delta_{X_j}$
  is an atomic measure on~$\cA$,
  while the uniform distribution on~$\cA$
  is absolutely continuous with respect to Lebesgue measure.
  The TV-norm of their difference does not vanish for $n\to\infty$:
  \begin{equation*}
    \TVnorm{\widehat\Prob_n-\Prob}
    =\sup_{A\in\Borel(\cA)}\abs{\widehat\Prob_n(A)-\Prob(A)}
    \ge1\text,
  \end{equation*}
  as the set $A:=\{X_n\mid n\in\N\}$ shows.
\end{Example}

We have to follow a different path.
Assume again $\cA\in\Borel(\R)$.
Let us abbreviate the difference of the cumulative distribution functions by
\begin{equation*}
  F_\Lambda(E)
  :=\spr{\ifu{\Ioc{-\infty}E}}{\widehat\Prob_\Lambda^\omega-\Prob_0}\textq,
  E\in\R\text,
\end{equation*}
and assume that~$g=f(\{0\},\argmt)$ has bounded variation,
or more specifically that it is monotone.
Then, we can perform the following partial integration
with Riemann-Stieltjes integrals:
\begin{align*}
  \abs{\spr g{\widehat\Prob_\Lambda^\omega-\Prob_0}_{\cA}}&
  =\Bigabs{\int_\cA g(E)\dd F_\Lambda(E)}
  =\Bigabs{-\int_\cA F_\Lambda(E)\dd g(E)}\\&
  \le\norm[\infty]{F_\Lambda}\int_\cA\d\abs g(E)
  =\TVnorm g\cdot\sup_{E\in\R}
    \abs{\spr{\ifu{\Ioc{-\infty}E}}{\widehat\Prob_\Lambda^\omega-\Prob_0}_\cA}
  \text.
\end{align*}
This calculations generalizes the theorem by Glivenko and Cantelli
to bounded monotone functions:  For $M>0$, we have
\begin{equation*}
  \sup\{\abs{\spr g{\widehat\Prob_{\Lambda_L}-\Prob_0}}
    \mid\text{$g\from\R\to\Icc{-M}M$ monotone}\}
  \xto{L\to\infty}0\text.
\end{equation*}

In order to deal with fields that are only \emph{almost} additive,
we have to treat patterns of all finite sizes and not only singletons.
Each pattern corresponds to a multivariate random variable.
This means that we require a multivariate version of Glivenko--Cantelli theory.
To formulate this, we introduce an multivariate version of the empirical measure: For given (large) set $\Lambda_j$, smaller set $\Lambda_L$, a grid $T_{j,L}$ and a coloring $\omega\in \cA^G$ we define
the empirical measure by
\[
 \widehat \Prob_{j,L}^\omega :\cB(\cA^{\Lambda_j}) \to [0,1],\qquad \widehat\Prob_{j,L}^\omega:=\frac{1}{|T_{j,L}|}\sum_{t\in T_{j,L}} \delta_{(\tau_t\omega)_{\Lambda_L}}.
\]
Here, the grid $T_{j,L}$ is a set of basepoints to (almost) cover the $\Lambda_j$ with translated versions of $\Lambda_L$ along $T_{j,L}$. An illustration of this it to be found in Figure \ref{fig:emp.measure}.
Let us emphasize that the illustration serves well in the $\Z^d$-case or in the ST-amenable case. 
However, for general amenable groups there is usually no grid $T_{j,L}$ for a (perfect) covering a set with \emph{one} set. In this case one can still use the above definition of the empirical measure, 
but, as in Section \ref{sec:general}, one needs to implement the technique of $\epsilon$-quasi tilings. Moreover, let use emphasize that counting patterns in the empirical measure along a grid is substantially different from counting patterns in the definition of frequencies in Section \ref{sec:patternsfrequencies}, see the definition of $\sharp_PP'$ in \eqref{eq:patterncount} and compare Figure \ref{fig:patterncount} with Figure \ref{fig:emp.measure}.
\begin{figure}
 
\begin{tikzpicture}[scale=0.6]
{
{\color{black}

{

\foreach \y in {0,1,...,8}{
\draw (1,\y) --(9,\y);
}
\foreach \x in {1,2,...,9}{
\draw (\x,0) --(\x,8);
}
\foreach \x in {2,4,6,8}{
\draw[dotted, thick] (\x+0.5,-0.5) -- (\x+0.5,8.5);
\draw[dotted, thick] (0.5,\x-0.5) -- (9.5,\x-0.5);
}

\filldraw[blue] (8,2) circle (4pt);
\filldraw[blue] (7,8) circle (4pt);
\filldraw[blue] (6,5) circle (4pt);
\filldraw[blue] (7,2) circle (4pt);
\filldraw[blue] (6,8) circle (4pt);

\filldraw[red] (1,1) circle (4pt);
\filldraw[green] (1,2) circle (4pt);
\filldraw[green] (1,3) circle (4pt);
\filldraw[red] (1,4) circle (4pt);
\filldraw[red] (1,5) circle (4pt);
\filldraw[red] (1,6) circle (4pt);
\filldraw[green] (1,7) circle (4pt);
\filldraw[green] (1,8) circle (4pt);
\filldraw[red] (2,1) circle (4pt);
\filldraw[green] (2,1) circle (4pt);
\filldraw[red] (2,2) circle (4pt);
\filldraw[green] (2,3) circle (4pt);
\filldraw[blue] (2,4) circle (4pt);
\filldraw[blue] (2,5) circle (4pt);
\filldraw[red] (2,6) circle (4pt);
\filldraw[red] (2,7) circle (4pt);
\filldraw[green] (2,8) circle (4pt);
\filldraw[green] (3,1) circle (4pt);
\filldraw[red] (3,2) circle (4pt);
\filldraw[red] (3,3) circle (4pt);
\filldraw[red] (3,4) circle (4pt);
\filldraw[green] (3,5) circle (4pt);
\filldraw[red] (3,6) circle (4pt);
\filldraw[blue] (3,7) circle (4pt);
\filldraw[green] (3,8) circle (4pt);
\filldraw[blue] (4,1) circle (4pt);
\filldraw[green] (4,2) circle (4pt);
\filldraw[green] (4,3) circle (4pt);
\filldraw[green] (4,4) circle (4pt);
\filldraw[blue] (4,5) circle (4pt);
\filldraw[blue] (4,6) circle (4pt);
\filldraw[blue] (4,7) circle (4pt);
\filldraw[blue] (4,8) circle (4pt);
\filldraw[red] (5,1) circle (4pt);
\filldraw[blue] (5,2) circle (4pt);
\filldraw[red] (5,3) circle (4pt);
\filldraw[green] (5,4) circle (4pt);
\filldraw[green] (5,5) circle (4pt);
\filldraw[blue] (5,6) circle (4pt);
\filldraw[blue] (5,7) circle (4pt);
\filldraw[blue] (5,8) circle (4pt);
\filldraw[red] (6,1) circle (4pt);
\filldraw[red] (6,2) circle (4pt);
\filldraw[green] (6,3) circle (4pt);
\filldraw[red] (6,4) circle (4pt);
\filldraw[green] (6,5) circle (4pt);
\filldraw[red] (6,6) circle (4pt);
\filldraw[red] (6,7) circle (4pt);
\filldraw[green] (6,8) circle (4pt);
\filldraw[green] (7,1) circle (4pt);
\filldraw[blue] (7,2) circle (4pt);
\filldraw[red] (7,3) circle (4pt);
\filldraw[red] (7,4) circle (4pt);
\filldraw[green] (7,5) circle (4pt);
\filldraw[red] (7,6) circle (4pt);
\filldraw[green] (7,7) circle (4pt);
\filldraw[green] (7,8) circle (4pt);
\filldraw[green] (8,1) circle (4pt);
\filldraw[blue] (8,2) circle (4pt);
\filldraw[green] (8,3) circle (4pt);
\filldraw[green] (8,4) circle (4pt);
\filldraw[blue] (8,5) circle (4pt);
\filldraw[blue] (8,6) circle (4pt);
\filldraw[red] (8,7) circle (4pt);
\filldraw[green] (8,8) circle (4pt);
\filldraw[red] (9,1) circle (4pt);
\filldraw[blue] (9,2) circle (4pt);
\filldraw[red] (9,3) circle (4pt);
\filldraw[green] (9,4) circle (4pt);
\filldraw[blue] (9,5) circle (4pt);
\filldraw[blue] (9,6) circle (4pt);
\filldraw[blue] (9,7) circle (4pt);
\filldraw[red] (9,8) circle (4pt);
\filldraw[blue] (1,0) circle (4pt);
\filldraw[blue] (2,0) circle (4pt);
\filldraw[red] (3,0) circle (4pt);
\filldraw[red] (4,0) circle (4pt);
\filldraw[green] (5,0) circle (4pt);
\filldraw[blue] (6,0) circle (4pt);
\filldraw[green] (7,0) circle (4pt);
\filldraw[red] (8,0) circle (4pt);
\filldraw[green] (9,0) circle (4pt);

}
}
}

\draw[rounded corners=0.1cm]
 (4.7,4.2)-- ++(0,1.1)-- ++(1.6,0)-- ++(0,-1.6)-- ++(-1.6,0)-- ++(0,0.5);
\draw[rounded corners=0.1cm]
 (0.7,2.2)-- ++(0,1.1)-- ++(1.6,0)-- ++(0,-1.6)-- ++(-1.6,0)-- ++(0,0.5);
\draw[rounded corners=0.1cm]
 (6.7,0.2)-- ++(0,1.1)-- ++(1.6,0)-- ++(0,-1.6)-- ++(-1.6,0)-- ++(0,0.5);

\end{tikzpicture}
\caption{The set $\Lambda_j$ is a square of side length~$9$ and $\Lambda_L$ a square of side length~$2$. 
$\Lambda_j$ is (almost) covered when translating $\Lambda_L$ along the positions of the dashed grid (given by $T_{j,L}$). 
The marked pattern $P$ is found at $3$ (of $16$ possible) positions along the grid. Thus, the empirical measure of this pattern is 
$\widehat \Prob_{j,L}^\omega(P)= \widehat \Prob_{9,2}^\omega(P)= \frac{3}{16}$.}
\label{fig:emp.measure}
\end{figure}
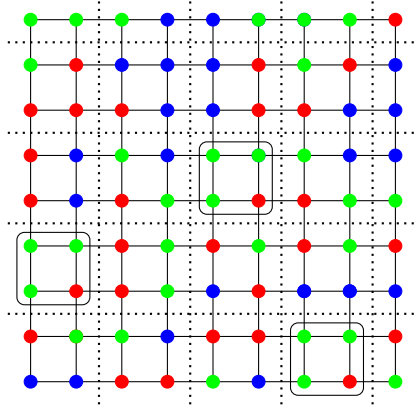

In order to apply the multivariate version of Glivenko--Cantelli,
we aim to integrate functions mapping from $\cA^\Lambda$ to~$\R$. 
Such a function $g\from\cA^\Lambda\to\R$ is called monotone,
if it is monotone in each coordinate.
Besides these generalizations due to higher dimensionality, 
there is another fundamental difference between univariate
and multivariate Glivenko--Cantelli theory:
While \Cref{thm:GC}
makes no assumptions on the distribution of the random variables,
the following example shows that we will have to impose some restrictions
on the joint distribution of the coordinates of the random vector.

\begin{Example}
  Let $X_j\sim\Normal(0,1)$, $j\in\N$,
  be i.\,i.\,d.\ standard normal random variables and $Y_j:=-X_j$.
  We consider the vectors $(X_j,Y_j)\in\R^2$.
  Let~$\widehat\Prob_n(\omega):=n^{-1}\sum_{j=1}^n\delta_{(X_j,Y_j)}$
  be the empirical distribution of $(X_j,Y_j)_{j=1}^n$ on~$\R^2$.
  The random test function
  \begin{equation*}
    g_\omega:=\ifu{\{(x,y)\in\R^2\mid x+y<0\}}
      +\ifu{\{(X_j(\omega),Y_j(\omega))\mid j\in\N\}}
    \from\R^2\to\Icc{-1}1
  \end{equation*}
  is monotone in each coordinate, and we have
  \begin{equation*}
    \sup\{\abs{\spr g{\widehat\Prob_n^\omega-\Prob}}\mid\text{$g\from\R^2\to\Icc{-1}1$
      monotone}\}
    \ge\abs{\spr{g_\omega}{\widehat\Prob_n^\omega-\Prob}}
    =1\text.
  \end{equation*}
\end{Example}

The problem arises because the set of discontinuities
of the monotone function has positive probability.
A correct generalization of \cref{thm:GC}
to the multivariate case is as follows.

\begin{Theorem}[DeHardt \cite{DeHardt-71}, Wright \cite{Wright-81}]\label{wright:LDP}
  Let
  \begin{itemize}
  \item $V_j$, $j\in\N$, be i.\,i.\,d.\ random variables with values in~$\R^k$
    and distribution~$\Prob$,
  \item $\widehat\Prob_n:=\frac1n\sum_{j=1}^n\delta_{V_j}$ for $n\in\N$
    the empirical distribution, and
  \item $M>0$ and $\cM:=\{g\from\R^k\to\Icc{-M}M\mid\text{$g$ monotone}\}$. 
  \end{itemize}
  Then, the following are equivalent.
  \begin{enumerate}[(i)]
  \item\label{GC-i}
    For all $J\subseteq\{1,\dotsc,k\}$, $J\ne\emptyset$,
    strictly monotone $g\from\R^J\to\R$, and $E\in\R$,
    the continuous part~$\Prob_c^J$ of the marginal~$\Prob^J$ of~$\Prob$
    satisfies
    \begin{equation*}
      \Prob_c^J\bigl(\partial g^{-1}\bigl(\Ioc{-\infty}E\bigr)\bigr)=0\text.
    \end{equation*}
  \item There exists an almost sure event~$\Omega_{\mathrm{unif}}$ on which
    \begin{equation*}
      \sup_{g\in\cM}\abs{\spr g{\widehat\Prob_n-\Prob}}\xto{n\to\infty}0\text.
    \end{equation*}
  \item For all $\kappa>0$, there are $a_\kappa,b_\kappa>0$ such that for all
    $n\in\N$, there is an event~$\Omega_{\kappa,n}$ with
    $\Prob(\Omega_{\kappa,n})\ge1-b_\kappa\exp(-a_\kappa n)$
    on which
    \begin{equation*}
      \sup_{g\in\cM}\abs{\spr g{\widehat\Prob_n-\Prob}}\le\kappa\text.
    \end{equation*}
  \end{enumerate}
\end{Theorem}

Here $\partial g^{-1}\bigl(\Ioc{-\infty}E\bigr)$
denotes the boundary of the sublevel set $g^{-1}\bigl(\Ioc{-\infty}E$.
Condition~\labelcref{GC-i} is trivial in the classical case $k=1$.
Also, in any dimension, each product measure~$\Prob$
satisfies condition~\labelcref{GC-i}.
In fact, the following theorem 
holds true.
\begin{Theorem}\label{thm:strictlymonotone}
  Let~$\Prob$ be a probability measure on~$\R^k$
  which is absolutely continuous with respect to a product measure
  $\Tensor_{j=1}^k\mu_j$ on~$\R^k$, where~$\mu_j$, $j\in\{1,\dotsc,k\}$
  are measures on~$\R$.
  Then, 
  condition \cref{wright:LDP}\labelcref{GC-i} is satisfied.
\end{Theorem}
See \cite[Theorem~5.5]{SchumacherSV-17} for a proof.

\subsection{Uniform limits for monotone fields}

The theorems which follow for the case of infinitely many colors~$\cA$
all have an additional assumption, namely the monotonicity in the random parameters.
This is a natural assumption, indeed: The IDS depends on the potential antitonely in our models.
Also, the IDS is monotone in site percolation.
\begin{Example}
  We revisit the Anderson model on a site percolation graph from
  \cref{sec:AndPerc}, this time with $\cA\subseteq\R$.
  Fix a bounded set $\cA_0\in\Borel(\R)$ for the values of the potential,
  and let $\cA:=\cA_0\union\{\alpha\}$ with $\alpha>2d+\sup\cA_0$.
  The value~$\alpha$
  of the potential is interpreted as a closed site in the percolation graph:
  \begin{equation*}
    \cV_\omega:=\{v\in G\mid V_\omega(v)\ne\alpha\}\text.
  \end{equation*}
  The edges of the percolation graph are as before
  \begin{equation*}
    \cE_\omega:=\{\{v,w\}\subseteq\cV_\omega\mid v\sim w\}\text.
  \end{equation*}
  The Hamiltonian $H_\omega\from\ell^2(G)\to\ell^2(G)$ is given by
  \begin{equation*}
    (H_\omega\varphi)(v):=
    \begin{cases}
      -\Laplace_\omega\varphi(v)+V_\omega\varphi(v)
        &\text{if $v\in\cV_\omega$, and}\\
      \alpha\varphi(v)
        &\text{if $v\in G\setminus\cV_\omega$.}
    \end{cases}
  \end{equation*}
  By the min-max principle,
  the eigenvalues do not decrease when we increase
  the potential at a site~$v\in G$.
  Particularly, when the potential reaches the value~$\alpha$
  and the site closes, the eigenfunction experiences de facto
  a Dirichlet boundary condition on that site,
  which also at most increases the kinetic energy.
  \par
  The eigenvalue counting functions count less eigenvalues
  below a given threshold, if the eigenvalues increase.
  Therefore, the eigenvalue counting functions decrease
  when the potential is raised.
  The same holds true for the limit, i.\,e.\ the IDS.
  This is the reason why the IDS in the quantum percolation model
  with random potential is antitone in the randomness.
\end{Example}

We first turn to the special case $G=\Z^d$.
It is physically most relevant and,
since~the group~$\Z^d$ satisfies the tiling property,
we do not need to resort to quasi tilings in this case.

\begin{Theorem}[\cite{SchumacherSV-17}]\label{thm:infiniteA_Zd}
  Let~$\cA\in\Borel(\R)$, $\Omega:=\cA^{\Z^d}$,
  and let $(\Omega,\Borel(\Omega),\Prob)$
  be a probability space such that~$\Prob$ satisfies
  \begin{itemize}
  \item $\Prob$ is translation invariant with respect to the $\Z^d$-action,
  \item for all $\Lambda\in\cF$, the marginal~$\Prob_\Lambda$
    is absolutely continuous with respect to a product measure on $\R^\Lambda$, and
  \item for a given $r\ge0$ and all $\Lambda_j\in\cF$, $j\in\N$,
    with $\min\{\dist(\Lambda_i,\Lambda_j)\mid i\ne j\}>r$,
    the random variables $\omega\mapsto\omega_{\Lambda_j}$, $j\in\N$,
    are independent.
  \end{itemize}
  Further, let $f\colon\cF\times\Omega\to\B$
  be a translation invariant, local, almost additive, monotone, bounded field.
  Then there exists an event~$\tilde\Omega$
  of full probability and a function $\fbar\in\B$
  such that for every $\omega\in\tilde\Omega$ we have
  \begin{equation}\label{eqmain}
    \lim_{j\to\infty}\biggnorm{
      \frac{f(\Lambda_j,\omega)}{\setsize{\Lambda_j}}-\fbar}
      =0\text,
  \end{equation}
  where $\Lambda_j:=\Ico0j\isect\Z^d$ for $j\in\N$.
  For an estimate on the speed of convergence,
  denote the bound of~$f$ by~$C_f$,
  the bound of the boundary term~$b$ of~$f$ by~$D_b$.
  Then, for every $\kappa>0$ and $L\in\N$, $L>2r$, there are $a,b>0$,
  depending on~$\kappa$, $L$, and~$C_f$,
  such that for all $j\in\N$, $j>2L$, there is an event~$\Omega_{\kappa,j}$
  with probability $\Prob(\Omega_{\kappa,j})\ge1-b\exp(-a\floor{j/L}^d)$,
  on which
  \begin{equation*}
      \Bignorm{\frac{f(\Lambda_j,\omega)}{\setsize{\Lambda_j}}- \fbar}
      \le2^{2d+1}\Bigl(\frac{(2C_f+D_b)L^d+D_br^d}{j-2L}
        +\frac{2(C_f+D_b)r^d+3D_br^d}{L-2r}\Bigr)
        +\kappa\text.
  \end{equation*}
  holds true.
\end{Theorem}

Of course, there is a version of \cref{thm:infiniteA_Zd} for amenable groups.
We follow the strategy in \cref{sec:general} and use quasi tilings
to deal with infinitely many colors on amenable groups.
This brings new challenges.
\Cref{wright:LDP} needs as input i.\,i.\,d.\ samples.
But quasi tilings are allowed to overlap (in a relatively small volume),
see \cref{def:quasitiling}.
This destroys the independence of eigenvalue counting functions associated to overlapping tiles.
One is tempted to excise the overlap from some of the tiles.
However, this would leave us with an independent but not identically distributed sample.
In this situation Glivenko--Cantelli-Theory is difficult to apply. 
The solution is to independently resample the portions of the quasi tiles
which overlap and to account for the error by a volume estimate.

The last result we present here treats fields with
infinitely many colors on amenable groups.

\begin{Theorem}[\cite{SchumacherSV-18}]\label{thm:main}
  Let~$G$ be a finitely generated amenable group
  with a F{\o}lner sequence~$(Q_j)_j$.
  Further, fix~$\cA\in\cB(\R)$, and let $(\Omega=\cA^G,\cB(\Omega),\Prob)$
  be a probability space such that~$\Prob$
  is translation invariant with respect to~$G$,
  has finite marginals with density w.\,r.\,t.\ a product measure,
  and independence at a distance.
  Further, let~$\cU$ be a set of translation invariant, local,
  almost additive, monotone, bounded fields $f\from\cF\times\Omega\to\B$
  with common bound, i.\,e.\ $C:=\sup\{C_f\mid f\in\cU\}<\infty$,
  common boundary term~$b\from\cF\to\R$ with bound~$D:=D_b$.
  \begin{enumerate}[(a), wide]
  \item Then, there exists an event $\tilde\Omega\in\cB(\Omega)$
    such that $\Prob(\tilde\Omega)=1$ and for any $f\in\cU$
    there exists a function~$\fbar\in\B$,
    which does not depend on the specific F{\o}lner sequence~$(Q_j)_j$, with
    \begin{equation*}
      \forall\omega\in\tilde\Omega\colon\quad
      \lim_{j\to\infty}
      \sup_{f\in\cU}\biggnorm{\frac{f(Q_j,\omega)}{\setsize{Q_j}}-\fbar}
      =0
      \text.
    \end{equation*}
  \item Furthermore, for each $\ve\in\Ioo0{1/10}$,
    there exist $j_0(\ve)\in\N$, independent of~$C$,
    such that for all $f\in\cU$, there are $a(\ve,C),b(\ve,C)>0$,
    such that for all $j\in\N$, $j\ge j_0(\ve)$,
    there is an event $\Omega_{j,\ve,C}\in\Borel(\Omega)$,
    with the properties
    \begin{equation*}
      \Prob(\Omega_{j,\ve,C})
      \ge1-b(\ve,C)
      \exp\bigl(-a(\ve,C)\setsize{Q_j}\bigr)
    \end{equation*}
    and
    \begin{align*}
      \biggnorm{\frac{f(Q_j,\omega)}{\setsize{Q_j}}-\fbar}&
      \le(37C_f+47\constb[\cU]+47)  \ve
      \qtext{ for all $\omega\in\Omega_{j,\ve,C}$ and all $f\in\cU$.}
    \end{align*}
  \end{enumerate}
\end{Theorem}

\section{Outlook}

We have presented a number of theorems concerning convergence in sup-norm
and other Banach-space norms of averaged almost additive fields.
More important than the convergence itself are the corresponding quantitative
error estimates.
They split into two parts of two different origins:
The geometric part and the probabilistic part.
The probabilistic error measures how far off certain empirical
measures are from their theoretical counterparts.
This difference can be estimated using large deviations techniques,
which is implicit in our use of the \cref{wright:LDP} of DeHardt and Wright.
An aspect which is not completely satisfactory is that we are not able
to specify the dependence of the positive coefficients~$a_\kappa$
and~$b_\kappa$ on the small parameter $\kappa>0$
and the dimension of the pattern $k\in\N$.
For this reason we are not able to choose the two lengths scales
which tend both to infinity as functions of each other.

Furthermore, the Theorems in \cref{sec:infiniteA}
assume certain monotonicity with respect to individual random parameters.
While this is sufficient for a wide variety of models in statistical physics,
there are examples, e.\,g.\ random hopping Hamiltonians,
which do not satisfy this assumption.
For this reason it is desirable to relax the monotonicity assumption,
or formulate alternative sufficient conditions.
These are aims which we will pursue in a forthcoming project.


\small 
\begin{longtabu} to\linewidth { r @{$\colon$} X[l] }
  \toprule
  \multicolumn{2}{c}{Table of Notation\label{tableofnotations}}\\
  \midrule
\endfirsthead
  \toprule
  \multicolumn{2}{c}{Table of Notation (continued)}\\
  \midrule
\endhead
  \midrule
  \multicolumn{2}{r}{\footnotesize continues on next page}\\
\endfoot
\bottomrule
\endlastfoot
$G$ & a finitely generated amenable group, or its Cayley graph.
\\ $g$                 & element of~$G$, when~$G$ is used as group
\\ $v,w$               & elements of~$G$, when~$G$ is used as Cayley graph
\\ $\tilde\tau_g$      & group action of the group~$G$ on its Cayley graph~$G$
\\ $U_g$               & unitary group action of~$G$ on $\ell^2(G)$
\\ $\Laplace$          & Laplace operator on $\ell^2(G)$, $-\Laplace\ge0$
\\ $\cA$               & set of colors
\\ $\Borel(\cA)$       & Borel sets on~$\cA$
\\ $\Prob_0$           & probability measure on $(\cA,\Borel(\cA))$
\\ $\Omega$            & set of colorings of~$G$
\\ $\omega$            & coloring
\\ $\Prob$             & probability measure on $(\Omega,\Borel(\Omega))$
\\ $V$                 & random potential
\\ $\tau_g$            & group action of~$G$ on~$\Omega$ and on patterns
\\ $H_\omega$          & random Schr\"odinger operator
\\ $\Sigma$            & almost sure spectrum of $H_\omega$
\\ $N$                 & integrated density of states (IDS)
\\ $\d N$              & density of states measure
\\ $\Lambda$           & (finite) subset of~$G$, often domain of pattern
\\ $\ell^2(\Lambda)$   & subspace of $\ell^2(G)$
\\ $\ifu\Lambda$       & projection $\ell^2(G)\to\ell^2(\Lambda)$
\\ $H_\omega^\Lambda$  & restriction of~$H_\omega$ to $\ell^2(\Lambda)$
\\ $\cF$               & set of finite subsets of~$G$
\\ $\delta_v$          & Kronecker delta on~$v$
\\ $n(\Lambda,\omega)$ & eigenvalue counting function, not normalized
\\ $\B$                & Banach space, often the right continuous $\R\to\R$-functions
\\ $\Lambda_L$         & cube of side length~$L$ or small F\o{}lner sequence
\\ $\Psite$            & Bernoulli measure for site percolation
\\ $\cV_\omega$        & vertices of percolation graph
\\ $\cE_\omega$        & edges of percolation graph
\\ $\Pedge$            & Bernoulli measure for edge percolation
\\ $e_j$               & $j$-th basis vector of $\Z^d$
\\ $\Laplace_\omega$   & Laplace operator on percolation graph
\\ $\alpha$            & eigenvalue of~$H_\omega$ on $G\setminus\cV_\omega$
\\ $\cA_0$             & set of values of the random potential on a percolation graph
\\ $\cV$               & set of visible points in~$\Z^d$
\\ $\Hvis$             & Schr\"odinger operator with~$\ifu\cV$ as potential
\\ $\tilde N$          & limiting function in \cref{thm:LMV-ecf}
\\ $P$                 & pattern, $P\from\Lambda\to\cA$
\\ $[P]_G$             & equivalence class of patterns with respect to~$G$
\\ $\omega_\Lambda$    & pattern induced by~$\omega$ on~$\Lambda$
\\ $P'|_\Lambda$       & pattern induced by $P'\in\cA^{\Lambda'}$ on~$\Lambda$
\\ $P'|_\cF$           & set of all patterns induced by~$P'$
\\ $\sharp_PP'$        & number of occurrences of~$P$ in~$P'$
\\ $\partial^r\Lambda$ & (two-sided) $r$-boundary of~$\Lambda$
\\ $\dist$             & graph distance on Cayley graph
\\ $(Q_j)_j$           & F\o{}lner sequence 
\\ $\nu_P$             & (asymptotic) frequency of~$P$
\\ $b$                 & boundary term
\\ $D_b$               & bound of~$b$
\\ $F$                 & field (without explicit dependence on~$\omega$)
\\ $\tilde F$          & pattern function of~$F$
\\ $C_F$               & bound of field~$F$
\\ $\Fbar$             & limit in \cref{thm:LenzMV}
\\ $\widehat\Prob_\Lambda^\omega$  & empirical distribution
\end{longtabu}


\def\polhk#1{\setbox0=\hbox{#1}{\ooalign{\hidewidth
  \lower1.5ex\hbox{`}\hidewidth\crcr\unhbox0}}}

\end{document}